\documentclass[10pt]{article}

\usepackage{amsfonts,amsmath,amssymb,amsthm}
\usepackage[tight]{subfigure}
\usepackage{kmac}
\usepackage[dvipdfm]{color}
\usepackage[final]{graphicx}
\usepackage{fullpage}

\usepackage{mathpple,courier}

\usepackage{calc}

\definecolor{purple}{rgb}{1.0,0.0,1.0}
\definecolor{orange}{rgb}{0.99,0.78,0.10}

\normalsize

\newcommand{\Basin}[1]{B(#1)}
\newcommand{\HH}{Hodgkin-Huxley}

\newcommand{\Imain}{14}
\newcommand{\iprod}[2]{\langle{#1},{#2}\rangle}

\newcommand{\lyap}{\varsub{\lambda}{max}}
\newcommand{\projss}{\varsub{\pi}{ss}}

\newcommand{\putgraph}[3]{\begin{center}\includegraphics[0in,0in][#1,#2]{#3}\end{center}}

\newcommand{\subgraph}[4]{\subfigure[#4]{\includegraphics[0.25in,0.25in][#1,#2]{#3}}}
\newcommand{\subgrapha}[5]{\subfigure[#4]{\includegraphics[bb=0in 0in #1 #2,scale=#3]{#5}}}

\title{\Large Entrainment and chaos in a pulse-driven {\HH} oscillator}

\author{Kevin K. Lin\\{\normalsize {\tt klin@cims.nyu.edu}}}

\date{January 18, 2006}

\begin{document}

\maketitle

\begin{abstract}
\noindent
The {\HH} model describes action potential generation in certain types
of neurons and is a standard model for conductance-based, excitable
cells.
Following the early work of Winfree and Best, this paper explores the
response of a spontaneously spiking {\HH} neuron model to a periodic
pulsatile drive.  The response as a function of drive period and
amplitude is systematically characterized.  A wide range of
qualitatively distinct responses are found, including entrainment to the
input pulse train
and persistent chaos.  These observations are consistent with a theory
of kicked oscillators developed by Qiudong Wang and Lai-Sang Young.  In
addition to general features predicted by Wang-Young theory, it is found
that most combinations of drive period and amplitude lead to entrainment
instead of chaos.  This preference for entrainment over chaos is
explained by the structure of the {\HH} phase resetting curve.
\end{abstract}



\section{Introduction}


The {\HH} model describes action potential generation in certain types
of neurons and is a standard model for conductance-based, excitable
cells {\cite{cronin,hh,hubel}}.  There is an extensive literature on the
response of the {\HH} model to different types of inputs
{\cite{aihara,alexander,glass,guttman,hayashi3,hayashi4,hayashi5,holden,matsumoto1,matsumoto2}},
and understanding how single neurons respond to external forcing
continues to be relevant for the study of information transmission in
neural systems {\cite{hunter-milton,mainen}}.
Because neurons typically communicate via pulsatile synaptic events, it
is natural to investigate the response of the {\HH} model to pulsatile
inputs.
Early studies by Best and Winfree {\cite{best,winfree}} examine the
response of a {\HH} model to periodic impulse trains, chracterizing in
detail the structure of phase singularities and the transition from
degree 1 to degree 0 phase resetting.  However, their work does not
systematically address the asymptotic dynamical behavior as a function
of drive period and amplitude.\footnote{\label{fn:takabe} Takabe,
Aihara, and Matsumoto {\cite{takabe}} appear to have carried out such a
systematic study.  But, I was only able to locate an abstract.}

This paper studies a spontaneously spiking ({\em i.e.} oscillatory)
{\HH} neuron model driven by periodic, pulsatile input of fixed
amplitude and period, and systematically classifies the response as a
function of drive period and amplitude.  It is found that:
\begin{enumerate}

\item In response to periodic pulsatile forcing of fixed amplitude $A$
  and period $T$, a spontaneously spiking {\HH} system can exhibit a
  wide range of distinct behaviors depending on $A$ and $T$:
  \begin{enumerate}

  \item\label{item:entrainment} {\bf Entrainment:} The driven system is
     stably periodic and its period is a rational multiple of the drive
     period $T$.

  \item\label{item:transient-chaos} {\bf Transient chaos:} The system
    experiences a transient period of exponential instability before
    entraining to the input.  This transient chaos is caused by a Smale
    horseshoe {\cite{gh}}.

  \item\label{item:chaos} {\bf Chaos:} The system becomes fully chaotic:
    it possesses a positive Lyapunov exponent and a mixing attractor
    (see {\cite{young}} for a review of these concepts).

  \end{enumerate}
  The response of the pulse-driven neuron is approximately
  $T_0$-periodic in the drive period $T$, where $T_0$ is the intrinsic
  period of the unforced {\HH} oscillator.  For example, if the
  pulse-driven oscillator is chaotic for some drive amplitude $A$ and
  drive period $T$, then it is likely to be chaotic when driven by a
  pulse train of amplitude $A$ with period near $T+T_0$.

\item The scenarios enumerated above are {\em prevalent} in the sense
  that they correspond to positive-area subsets of the drive
  period-drive amplitude space.  Prevalence, together with the
  approximate periodicity stated above, imply that each scenario occurs
  with positive ``probability.''  (See the discussion of
  Fig.~\ref{fig:lyaps} in {\S\ref{sec:results}} for the precise meaning
  of probability in this context.)  The range of responses and their
  prevalence are consistent with a theory of nonlinear oscillators
  developed recently by Qiudong Wang and Lai-Sang Young
  {\cite{wy1,wy2,wy3,wy4}}.

\item While chaotic behavior is readily observable, most combinations of
  drive period and drive amplitude lead to entrainment instead of chaos.
  This preference for entrainment can be explained by the structure of
  the phase resetting curve (see {\S\ref{sec:theory}}) of the {\HH}
  system.

\end{enumerate}


This paper relies heavily on numerical computation and the conceptual
framework provided by Winfree's theory of biological rhythms
{\cite{winfree}} and the work of Q. Wang and L.-S. Young on nonlinear
oscillators {\cite{wy2,wy3}}.  Phase resetting curves, introduced by
Winfree, play a particularly important role here.  The phase resetting
curve of a nonlinear oscillator is an interval map which captures the
asymptotic response of a nonlinear oscillator to a single, pulsatile
perturbation.  Because they are $1$-dimensional objects, phase resetting
curves are often easier to understand than the nonlinear oscillators
they represent.  They are frequently used to infer stable dynamical
behavior like phase locking.  Wang-Young theory provides a mathematical
framework for using phase resetting curves to infer the existence and
prevalence of chaotic behavior.  Rather than numerically verify the
hypotheses of their theorems, we have opted to examine the consequences
of the theory directly, relying on a combination of numerical simulation
and geometric reasoning to characterize the specific response of the
{\HH} model to a periodic pulsatile drive.

For the sake of clarity, parameters are selected to ensure that the
{\HH} system possesses a unique limit cycle and no other attracting
invariant set.  This corresponds to a repeatedly spiking neuron with an
unstable rest state.
While the scenarios stated above should still hold when the limit cycle
coexists with other stable invariant sets, this choice simplifies the
interpretation of numerical simulations.  Otherwise, a trajectory may
jump out of the basin of the limit cycle, which obscures the mechanism
described by Wang-Young theory and which Winfree and Best have already
investigated thoroughly {\cite{best,winfree}}.

The rest of this paper is organized as follows: Section
{\ref{sec:unforced-hh}} briefly reviews the unforced {\HH} equations and
its properties.  Main numerical results are summarized in
{\S\ref{sec:results}} and discussed in {\S\ref{sec:theory}}.  Section
{\ref{sec:more-results}} discusses further numerical results, addressing
some issues raised in Sections {\ref{sec:results}} and
{\ref{sec:theory}}.  Section {\ref{sec:future}} discusses possible
extensions and generalizations.



\section{Brief review of the {\HH} model}
\label{sec:unforced-hh}



\begin{figure}
\putgraph{4in}{2.5in}{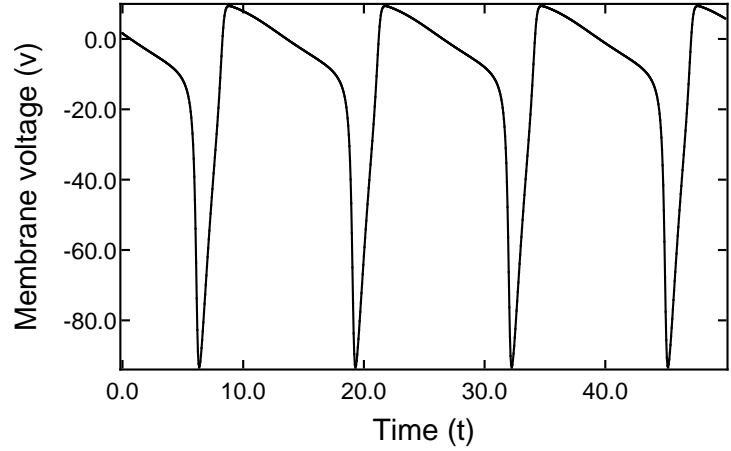}
\caption{The time course for the {\HH} equations at the parameter values
  (\ref{eqn:hh-params}).  The rapid ``spike'' followed by a long
  ``recovery'' period is typical of the {\HH} equations.}
\label{fig:hh-time-course}
\end{figure}

The {\HH} equations are a system of nonlinear ordinary differential
equations\footnote{This paper does not treat the {\HH} PDEs: spatial
dependence is not relevant here.}  which describe the way neurons
generate spatially and temporally localized electrical pulses
{\cite{cronin,hh,hubel}}.  These electrical pulses, called {\em action
potentials}, are the primary way in which neurons transmit information.
Action potentials are triggered by sufficiently large membrane voltages,
which can be set up by the influx of ions into the cell.
A neuron is said to {\em fire} or {\em spike} when it generates an
action potential (Fig. {\ref{fig:hh-time-course}}).  The {\HH} model
describes action potential generation in terms of the membrane voltage
and dimensionless {\em gating variables} which quantify the effective
permeability (or {\em conductance}) of the membrane for various types of
ions.


The original {\HH} equations model action potential generation in the
squid giant axon.  This giant axon contains two types of membrane ion
channels.  One type of channel is specific to potassium ions, the other
to sodium ions.  The state variables of the model are the membrane
voltage $v$, the activation $n$ of the potassium channels, and the
activation $m$ and inactivation $h$ of the sodium channels.  The
equations are {\cite{hh}}
\begin{equation}
\label{eqn:hh}
\begin{array}{ccl}
\dot{v} &=& C^{-1}\brac{-I - \varsub{\bar{g}}{K}n^4(v-\varsub{v}{K}) -
\varsub{\bar{g}}{Na}m^3h(v-\varsub{v}{Na}) -
\varsub{\bar{g}}{leak}(v-\varsub{v}{leak})}\\
\dot{m} &=& \varsub{\alpha}{m}(v)(1-m) - \varsub{\beta}{m}(v)m\\
\dot{n} &=& \varsub{\alpha}{n}(v)(1-n) - \varsub{\beta}{n}(v)n\\
\dot{h} &=& \varsub{\alpha}{h}(v)(1-h) - \varsub{\beta}{h}(v)h\\
\end{array}
\end{equation}
where
\begin{equation}
\begin{array}{ll}
\varsub{\alpha}{m}(v) = \Psi\of{\frac{v+25}{10}},
&\varsub{\beta}{m}(v)  = 4\exp\of{v/18},\\
\varsub{\alpha}{n}(v) = 0.1\Psi\of{\frac{v+10}{10}},
&\varsub{\beta}{n}(v)  = 0.125\exp\of{v/80},\\
\varsub{\alpha}{h}(v) = 0.07\exp\of{v/20},
&\varsub{\beta}{h}(v)  = \frac{1}{1+\exp\of{\frac{v+30}{10}}},\\
\Psi(v) = \frac{v}{\exp(v)-1}.&\\
\end{array}
\end{equation}
Each ion channel consists of independent, identical subunits which must
all open to allow ions to pass through.  The gating variables $m$, $n$,
and $h$ take value in $(0,1)$ and represent the fraction of subunits
which are open.  The term $n^4$ enters into the potassium conductance
because potassium channels consist of 4 identical subunits; analogous
structures account for the $m^3h$ term in the sodium conductance
{\cite{cronin}}.
The gating variable equations are master equations for continuous-time
Markov chains with voltage-dependent transition rates $\alpha$ and
$\beta$; the Markov chains describe the opening and closing of the
corresponding channel subunits.  The $\dot{v}$ equation is Kirchoff's
current law.  Action potentials are {\em downward} voltage spikes and a
positive $I$ corresponds to an {\em inflow} of positively-charged ions.
The voltage convention here is that of {\cite{hh}} and opposite
contemporary usage: the membrane voltage $v$ is defined by
\begin{displaymath}
v =\mbox{voltage outside} -\mbox{voltage inside}.
\end{displaymath}

Action potentials are generally initiated by perturbations to the
membrane voltage.  Such perturbations may be caused, for instance, by
the flow of ions across the cell membrane.  Because neurons transmit
signals through spatially and temporally localized pulses, it is natural
to model stimuli as impulses {\cite{rieke}}.
The simplest type of repetitive, pulsatile stimulus to a neuron is a
periodic impulse train.  This means replacing the $\dot{v}$ equation
above by
\begin{align}
\label{eqn:pulse-driven-hh}
\dot{v} =& C^{-1}\brac{-I - \varsub{\bar{g}}{K}n^4(v-\varsub{v}{K}) -
\varsub{\bar{g}}{Na}m^3h(v-\varsub{v}{Na}) -
\varsub{\bar{g}}{leak}(v-\varsub{v}{leak})}\\\nonumber
&+ A\sum_{k\in\Z}{G(t-kT)},
\end{align}
where $G$ is a ``bump'' function such that $\int{G(t)\ dt}=1$.  For
simplicity, most of this paper uses the choice $G(t)=\delta(t)$; Section
{\ref{sec:other-params}} discusses the response of the {\HH} system to a
pulsatile drive with
\begin{equation}
G(t) =\left\{\begin{array}{ll}
1/t_0,&0\leq t\leq t_0\\
0,&\mbox{otherwise}\\
\end{array}\right..
\end{equation}
Mathematically, one can also choose to perturb the gating variables, but
such perturbations are not entirely natural and are not considered here.

This paper uses the original {\HH} parameters {\cite{hh}}:
\begin{equation}
\label{eqn:hh-params}
\begin{array}{ll}
\varsub{v}{Na} = -115\mbox{ mV},
&\varsub{\bar{g}}{Na} = 120\mbox{ m$\Omega^{-1}$/cm}^2,\\
\varsub{v}{K} = +12\mbox{ mV},
&\varsub{\bar{g}}{K} = 36\mbox{ m$\Omega^{-1}$/cm}^2,\\
\varsub{v}{leak} = -10.613\mbox{ mV},
&\varsub{\bar{g}}{leak} = 0.3\mbox{ m$\Omega^{-1}$/cm}^2,\\
C = 1\mbox{ $\mu$F/cm}^2.&\\
\end{array}
\end{equation}
Time is measured in milliseconds and current density in $\mu$A/cm$^2$.

\begin{figure}
\begin{center}
\putgraph{5in}{5in}{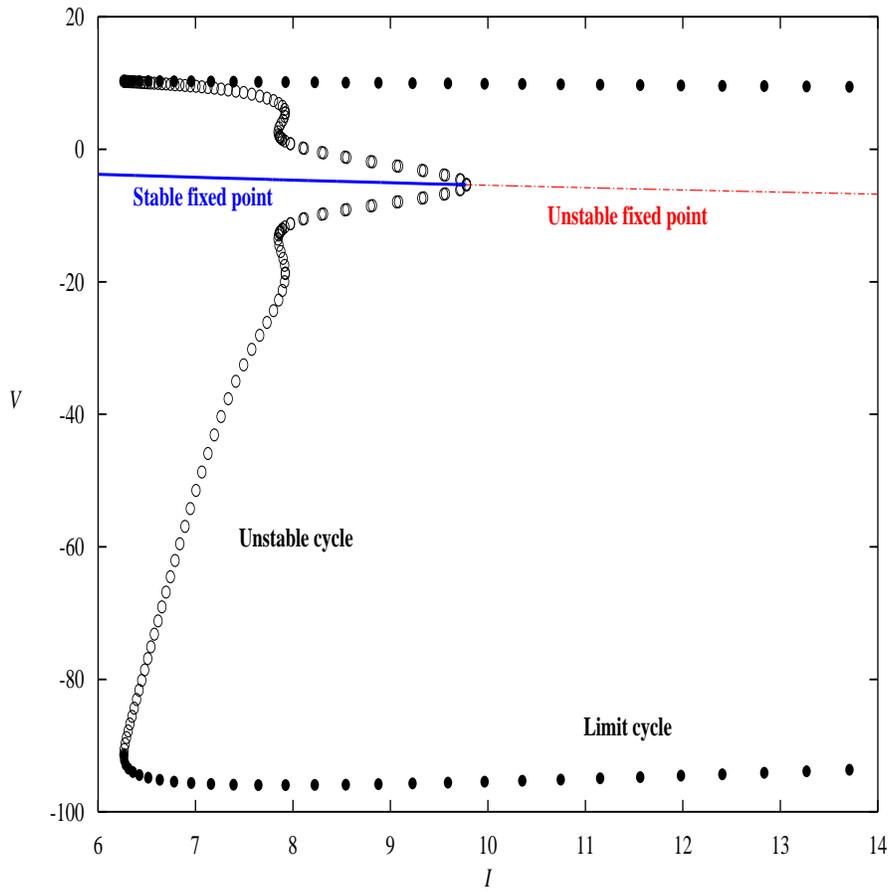}
\end{center}
\caption{The bifurcation diagram for the {\HH} equations as the injected
  current $I$ is varied.  The line in the middle marks the $v$
  coordinate of the rest state.  The solid blue part is stable while the
  dashed red part is unstable.  Solid black dots near the top and the
  bottom of the figure are the maximum and minimum $v$ values of limit
  cycles.  Empty black circles are the maximum and minimum $v$ values of
  unstable cycles.  The fixed point undergoes a subcritical Hopf
  bifurcation as $I$ increases.  This diagram is computed using XPPAUT
  {\cite{ermentrout}}.}
\label{fig:bifurcations}
\end{figure}
Figure {\ref{fig:bifurcations}} shows a bifurcation diagram for the
unforced {\HH} equations.  When $I=0$, the neuron maintains a stable
rest state, corresponding to the branch of stable fixed points on the
left of the diagram.  A sufficiently large value of $I$ causes a neuron
to fire repeatedly, which corresponds to the creation of a limit cycle
through a saddle-node bifurcation of periodic orbits.  Further
increasing $I$ destablizes the rest state through a subcritical Hopf
bifurcation.

In this paper, the injected current is always set to a value near
$I\approx\Imain$, corresponding to a steady ionic current which
destabilizes the rest state.  The phenomena studied here are insensitive
to the precise value of $I$ as long as it ensures the existence of a
stable limit cycle and an unstable fixed point.
As explained in the Introduction, these properties simplify the
interpretation of numerical simulations.  For this choice of $I$, the
Jacobian of the {\HH} vector field (Eq. {\ref{eqn:hh}}) at the unstable
fixed point has two real eigenvalues $\set{-4.97815, -0.146991}$ in the
left half plane, and a complex conjugate pair $0.0763367\pm 0.61866i$ in
the right half plane.  The fixed point thus has $2$-dimensional stable
and unstable manifolds.  The Lyapunov exponents associated with the
limit cycle are $0$, $\approx -0.20$, $\approx -2.0$, and $\approx
-8.3$.  Its period is $T_0\approx 12.944$.











\section{Main numerical results}
\label{sec:results}

\begin{figure}
\begin{center}
\includegraphics[bb=0in 0in 3in 5in]{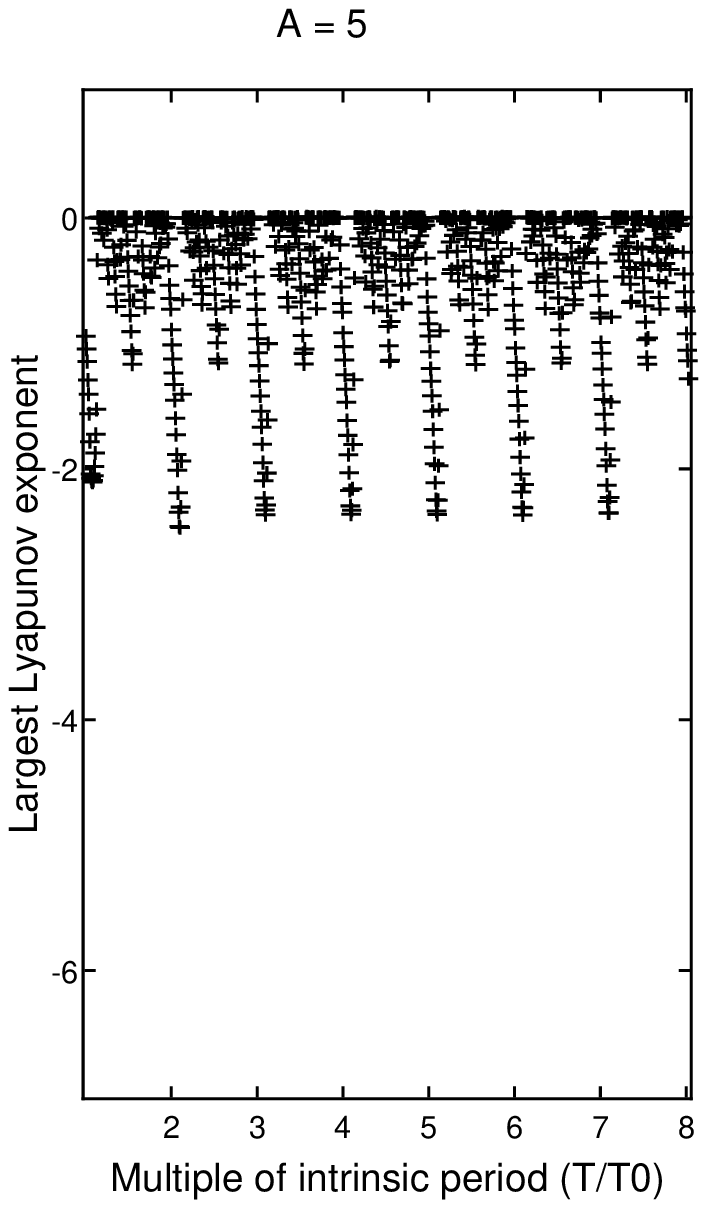}\includegraphics[bb=0in 0in
  3in 5in]{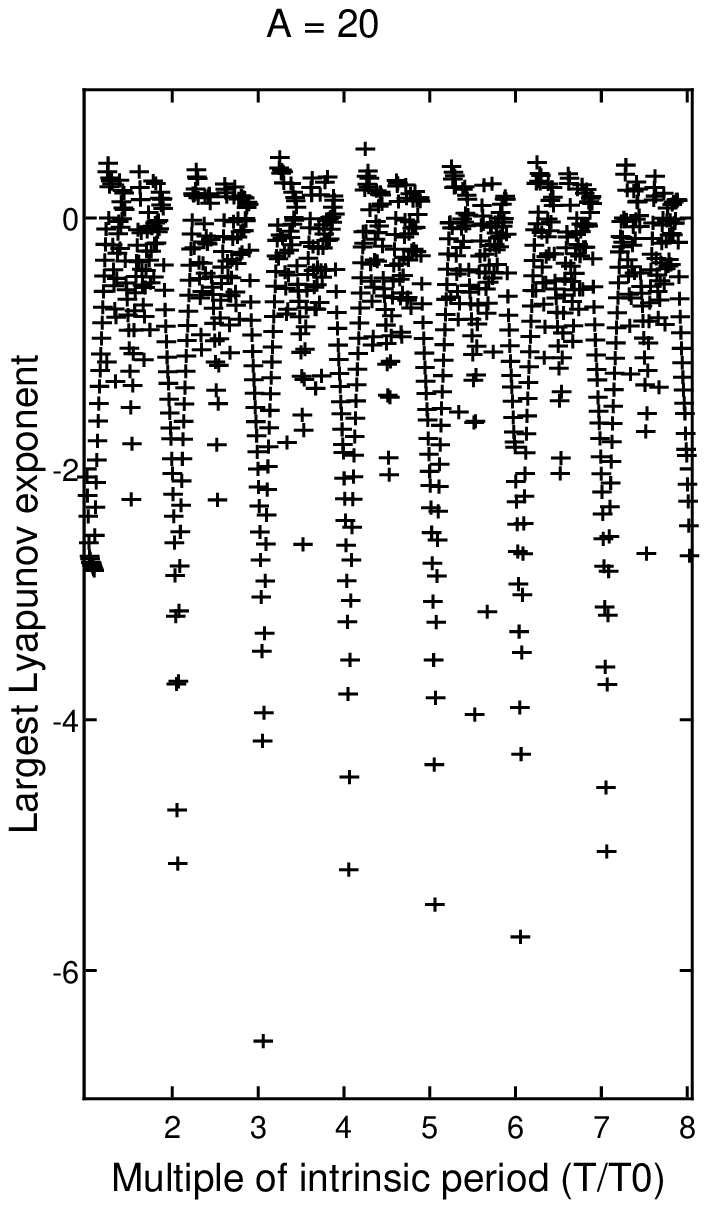}
\end{center}
\caption{Asymptotic properties of the pulse-driven flow are described by
the dynamics of the time-$T$ map $F_T$ (see Eq.~\ref{eqn:F_T}) and its
largest Lyapunov exponent $\lyap$.  Entrainment corresponds to $\lyap <
0$, and chaos corresponds to $\lyap > 0$.  This figure shows $\lyap$ as
a function of the drive period $T$, with $T$ ranging from $T_0\approx
13$ (the intrinsic period of the unforced {\HH} system; see
{\S\ref{sec:unforced-hh}}) to $8\cdot T_0\approx 101$.  {\bf Left:} Kick
amplitude is $A=5$.  {\bf Right:} Kick amplitude is $A=40$.  Note (i)
$\lyap(T+T_0)\approx\lyap(T)$; (ii) presence of both positive and
negative exponents for strong kicks (right), and only zero and negative
exponents for weak kicks (left); and (iii) the presence of more negative
exponents than positive ones.  See {\S\ref{sec:theory}} for a
discussion.  Lyapunov exponents are estimated by iterating $F_T$ for
1000 steps and tracking the rate of growth of a tangent vector.}
\label{fig:lyaps}
\end{figure}

\begin{figure}
\putgraph{5in}{3.5in}{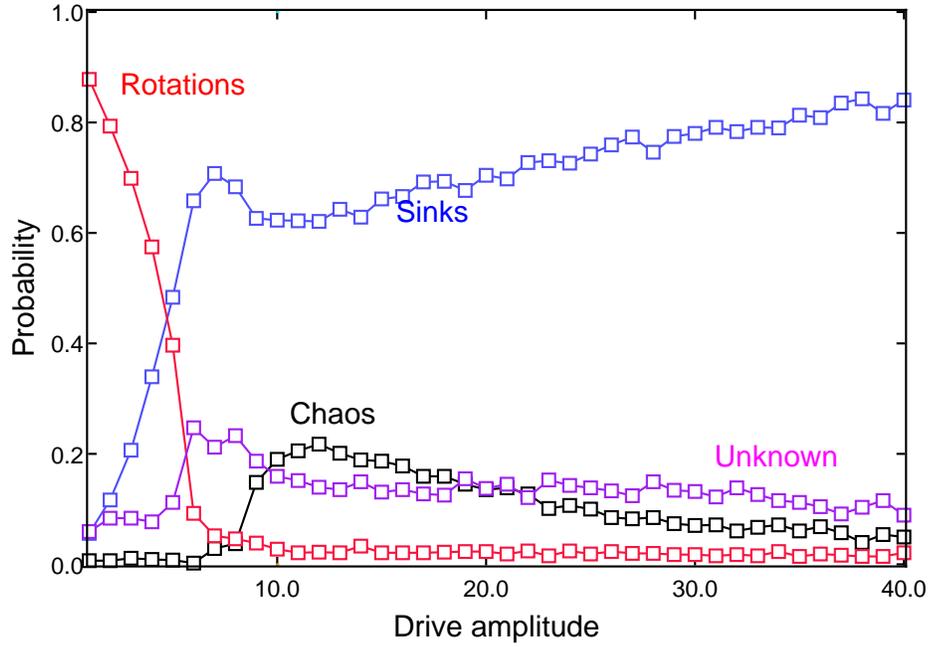}
\caption{The probability of different response types, as a function of
the drive amplitude $A$.
For each drive amplitude $A$, the fraction of drive periods $T\in[T_0,
8T_0]$ for which $\lyap\of{F_T} > 0$, {\em etc.}, is computed by
sampling from a uniform grid in $[T_0, 8T_0]$.  It is natural to equate
these fractions with probabilities because the Lyapunov exponents are
roughly periodic functions of $T$ (and become more so as $T\to\infty$),
as shown in Figure {\ref{fig:lyaps}} and explained in
{\S\ref{sec:theory}}.  {\em Empirical definitions:} Let $\hat\lambda$
denote the estimated Lyapunov exponent and $\eps$ the estimated standard
error.  Then ``chaos'' is defined as $\hat\lambda > 3\eps$,
``entrainment'' $\hat\lambda < -3\eps$, and ``rotation''
$\abs{\hat\lambda} <\eps/3$.}
\label{fig:lyaps-probs}
\end{figure}

Lyapunov exponents provide a convenient way to characterize the
asymptotic dynamics of Eq.~\ref{eqn:pulse-driven-hh}.  Let
$\phi_t:\R^4\to\R^4$ denote the flow map generated by the unforced {\HH}
equations, $T$ the drive period, and $A$ the drive amplitude.  The
Poincar\'e map
\begin{equation}
\label{eqn:F_T}
F_T(v,m,n,h) =\phi_T(v+A,m,n,h)
\end{equation}
takes a {\HH} state vector $(v,m,n,h)$, applies a pulse of amplitude $A$
to the membrane voltage, then evolves it for time $T$.  Iterating the
map $F_T$ thus gives a stroboscopic record of the state of our
pulse-driven {\HH} system before the arrival of each pulse.
The long-term dynamical behavior of the pulse-driven {\HH} oscillator
can be deduced from the asymptotic dynamics of $F_T$, which is
characterized by its (largest) Lyapunov exponent $\lyap$ {\cite{gh}}:
\begin{displaymath}
\begin{array}{lll}
\lyap < 0 &\Leftrightarrow\mbox{ $F_T$ has sinks}
&\Leftrightarrow\mbox{ kicked flow is entrained to input}\\
\lyap = 0 &\Leftrightarrow\mbox{ $F_T$
is quasiperiodic}
&\Leftrightarrow\mbox{ kicked flow drifts relative to input}\\
\lyap > 0 &\Leftrightarrow\mbox{ $F_T$ chaotic}
&\Leftrightarrow\mbox{ kicked flow is chaotic}\\
\end{array}
\end{displaymath}
(Sinks refer to stable fixed points and stable periodic orbits.)  Note
that of the scenarios given in the Introduction, transient chaos alone
does not appear in this list: Lyapunov exponents, being long-time,
average quantities, cannot detect transient chaos.

Figure {\ref{fig:lyaps}} shows the Lyapunov exponents of $F_T$ as a
function of $T/T_0$, where $T_0$ is the period of the {\HH} limit cycle.
Different colors correspond to different values of $A$.  The periodicity
of the response as a function of $T$ is apparent.  Because the response
type as a function of $T$ is approximately identical over each period
$[nT_0,(n+1)T_0]$, it makes sense to speak of the probability that a
randomly chosen drive period $T$ will elicit a particular asymptotic
behavior, for example chaos.  More precisely, periodicity ensures that
the fraction $p_n$ of drive periods $T$ in $[0,nT_0]$ for which $\lyap >
0$ converges to a well-defined limit as $n\to\infty$.  Similar
statements hold for $\lyap < 0$ and $\lyap=0$.

Figure {\ref{fig:lyaps-probs}} shows these probabilities as functions of
$A$.  At $A=10$, the probability of obtaining a positive exponent is
roughly 20\% and the probability of obtaining a negative exponent is
roughly 70\%.  Thus, if one were to pick $T$ randomly out of an interval
$[NT_0, (N+1)T_0]$ for large, fixed integer $N$, the probability that
$\lyap\of{F_T} > 0$ is about 20\%.  Figure {\ref{fig:lyaps-probs}} shows
that when $A$ is small, the most likely type of behavior is
rotation-like behavior.  This possibility becomes less likely as $A$
increases.  At the same time, sinks and chaos both become more likely,
with sinks dominating the scene.  One feature of Figure
{\ref{fig:lyaps-probs}} specific to the pulse-driven {\HH} flow is that
when $A$ is large, the system prefers entrainment over chaos in the
sense that entrainment has higher probability.  This preference is more
pronounced as $A$ increases.
Note that in computing Lyapunov exponents numerically, we only have
access to finite time information.  In principle, this means it is
virtually impossible to distinguish persistent chaotic behavior from
transient chaos caused by a ``large'' horseshoe (but see
{\S\ref{sec:horseshoe}}).

In all numerical simulations shown in this paper, Eq.~\ref{eqn:hh} is
integrated using an adaptive integrator of Runge-Kutta-Fehlberg type,
with an error tolerance of $10^{-6}$ (in the sup norm) {\cite{recipes}}.
The largest Lyapunov exponent $\lyap$ of $F_T$ is computed in a
straightforward manner, by choosing a nonzero unit vector $w\in\R^4$ and
estimating the rate of growth of $\norm{(DF_T)^nw}$.  The matrix-vector
product $(DF_T)^nw$ is easily computed via the variational equations
$\dot{x} = H(x),\dot{w} = DH(x)\cdot w$ for the {\HH} vector field $H$
($DH$ is the Jacobian matrix of $H$; see {\cite{geist}} for details).


\section{Discussion}
\label{sec:theory}

\subsection{Response to a single pulse: phase resetting curves}

This section reviews phase resetting curves.  See Winfree
{\cite{winfree}}, Glass and Mackey {\cite{glass}}, and Brown, Moehlis,
and Holmes {\cite{brown}} for more details and applications, and
Guckenheimer and Holmes {\cite{gh}} for background information on
dynamical systems theory.  See
{\cite{ermentrout,ek,ermentrout1,williams-bowtell}} for further
discussions of phase resetting curves.

Let $\phi_t:\R^n\to\R^n$ be a flow generated by a smooth vector field
with a hyperbolic limit cycle $\gamma$.  Such a limit cycle represents a
stable nonlinear oscillator.  The basin of attraction of $\gamma$ is
denoted $\Basin{\gamma}$.  The hyperbolicity of $\gamma$ guarantees that
points in $\Basin{\gamma}$ converge to $\gamma$ exponentially fast.  (It
is convenient to use $\gamma$ to refer to both the trajectory
$\gamma:\R\to\R^n$ and the invariant point set it defines.)  An
impulsive perturbation (``kick'') to the nonlinear oscillator can be
defined by specifying a kick amplitude $A$ and a kick direction
$\hat{K}:\R^n\to\R^n$ and defining a family of {\em kick maps}
\begin{equation}
K_A(x) = x + A\cdot\hat{K}(x),
\end{equation}
so that kicks send each point $x\in\R^n$ to $K_A(x)$.  For what follows,
$K_A$ should be smooth and satisfy
$K_A(\Basin{\gamma})\subset\Basin{\gamma}$.

The {\HH} system with the value of $I$ given in
{\S\ref{sec:unforced-hh}} is a nonlinear oscillator whose basin
$\Basin{\gamma}$ is an open subset of $\R^4$.  The kick map
corresponding to an instantaneous voltage spike is simply $K_A(v,m,n,h)
= (v+A, m, n, h)$.  As in {\S\ref{sec:results}}, it is convenient to
introduce the time-$T$ map
\begin{equation}
F_T =\phi_T\circ K_A,
\end{equation}
where $\circ$ denotes function composition.
Iterating $F_T$ gives a stroboscopic record of the system state before
the arrival of each kick, and thus describes the long-time dynamics of
the flow $\phi_t$ under repeated, $T$-periodic kicks.

Because the phase dimension $n$ may be large, the dynamics of
$F_T:\R^n\to\R^n$ may be difficult to analyze.  Winfree observed that
every point near the limit cycle $\gamma$ must converge to $\gamma$ as
$t\to\infty$, so the flow near $\gamma$ is dominated by the rotational
motion along $\gamma$.  Thus, one can reduce the dimension of the phase
space from $n$ to $1$, at least heuristically.  To do this, first define
the phase function $\theta:\gamma\to[0,T_0)$ by fixing a reference point
$x_0\in\gamma$ and requiring that $\phi_{\theta(x)}(x_0) = x$ for all
$x\in\gamma$.  By construction, $\theta$ satisfies
$\frac{d}{dt}\of{\theta(\gamma(t))} = 1, 0\leq t<T_0$.  The function
$\theta$ can be extended to a function $\theta:\Basin{\gamma}\to[0,T_0)$
by projecting along strong-stable manifolds\footnote{The
{\em strong-stable manifold} $\varsub{W}{ss}(x)$ of $x\in\gamma$ is the
set
\begin{equation}
\varsub{W}{ss}(x) =\set{y\in\R^n :
\lim_{n\in\Z,n\to+\infty}\phi_{nT_0}(y)\to x}.
\end{equation}
When the vector field generating $\phi_t$ is smooth, the strong-stable
manifolds are (locally) smooth submanifolds of $\R^n$.  The
strong-stable linear subspace $\varsub{E}{ss}(x)$ is the tangent space
of $\varsub{W}{ss}(x)$ at $x$.  See {\cite{guckenheimer,gh}}.}: if $y$
is a point in the basin of $\gamma$ then $\theta(y)$ is defined to be
$\theta(x)$, where $x$ is the unique point such that
$y\in\varsub{W}{ss}(x).$ This definition of phase preserves the property
that $\frac{d}{dt}\of{\theta(\phi_t(x))} = 1$.

Consider the limit {\cite{guckenheimer}}
\begin{equation}
\label{eqn:phase-resetting-def}
\bar{F}_T = \lim_{n\to+\infty}{F_{T+nT_0}}.
\end{equation}
The map $\bar{F}_T$ is well-defined on the basin of $\gamma$ and
retracts the basin onto $\gamma$, {\em i.e.} $\bar{F}_T(x)\in\gamma$ for
all $x\in\Basin{\gamma}$.  Thus, $\bar{F}_T$ induces an interval map
$f_T:[0,T_0)\to[0,T_0)$ which, given the current phase of the system,
yields the new phase after kicking and evolving the system for time $T$.
That is, $f_T(\theta(x)) =\theta(F_T(x))$ for all $x\in\Basin{\gamma}$.


The map $f_T$ is the {\em phase resetting curve}\footnote{Wang and Young
refer to phase resetting curves as {\em singular limits}.  Phase
resetting curves are also sometimes called {\em phase transition curves}
{\cite{glass}}.}, or more precisely the finite phase resetting curve
({\em infinitesimal phase resetting curves} {\cite{brown,ermentrout}}
are not needed here).  By construction, it has the property that
\begin{equation}
\label{eqn:prc-periodicity}
f_{T+\delta}(t) = f_T(t)+\delta\mbox{ (mod $T_0$)}.
\end{equation}
Thus, the family of maps $\set{f_T}$ is periodic in $T$.


\begin{figure}
\begin{center}
\subgraph{3in}{3in}{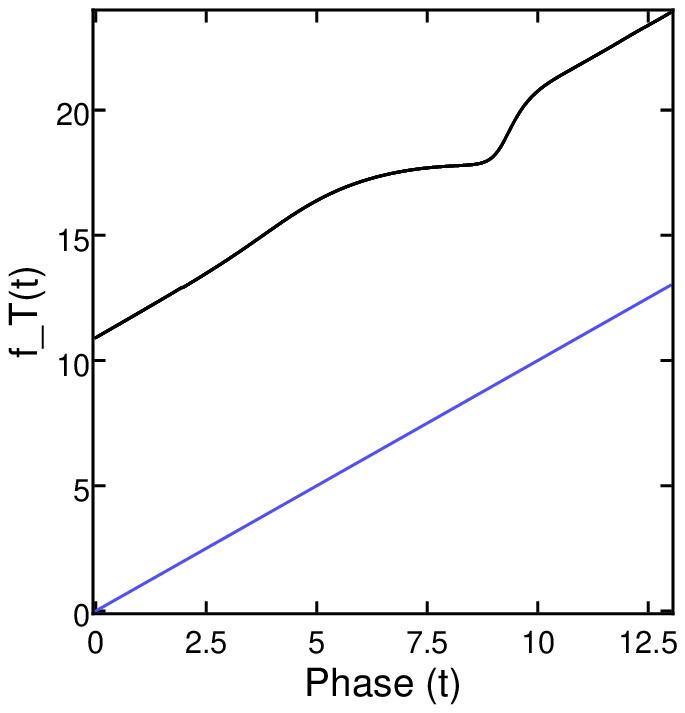}{$A=5$}\qquad
\subgraph{3in}{3in}{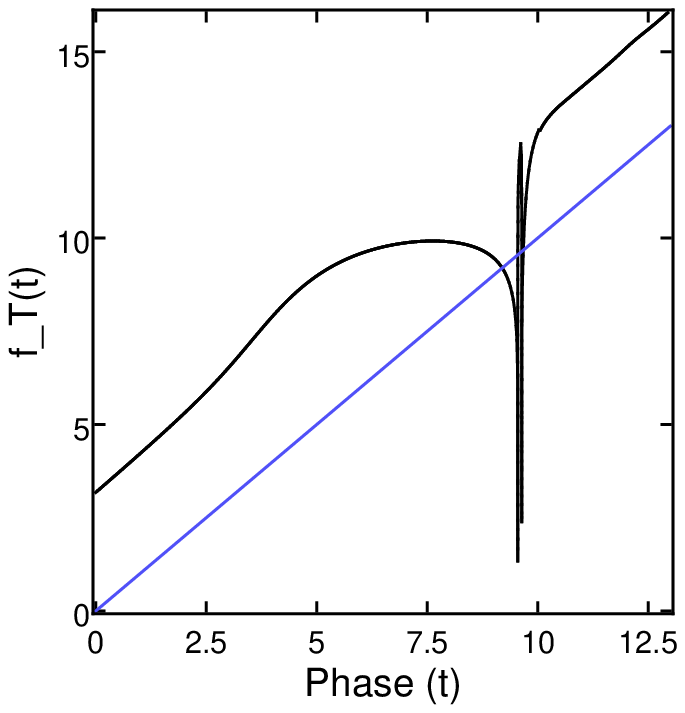}{$A=10$}\\
\subgraph{3in}{3in}{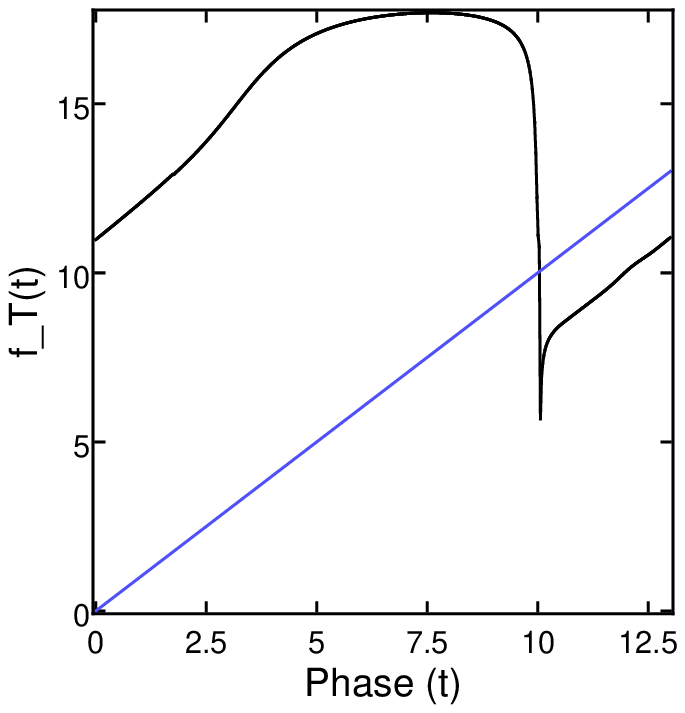}{$A=20$}
\end{center}
\caption{The graph of the lift $\tilde{f}_T$ of $f_T$, {\em i.e.} the
unique continuous map $\R\to\R$ such that $\tilde{f}_T = f_T$ on
$[0,T_0)$ and $\tilde{f}_T(t+T_0) =\tilde{f}_T(t)\mbox{ (mod $T_0$)}$,
for the pulse-driven {\HH} equations.  Drive amplitudes are (a) $A=5$,
(b) $A=10$, and (c) $A=20$.  The precise value of the drive period $T$
is not so important; varying $T$ shifts the graph vertically
(Eq. {\ref{eqn:prc-periodicity}}).  Note that $f_T$ has winding number 1
for $A=5$ and $A=10$, and has winding number 0 for $A=20$.  The
numerical data suggests that the degree changes around $A\approx
13.589$; the precise geometric mechanism is not clear.  Note that phase
ranges from $0$ to the intrinsic period $T_0\approx 12.9$ of the {\HH}
limit cycle.}
\label{fig:prc1}
\end{figure}

\begin{figure}
\putgraph{3in}{3in}{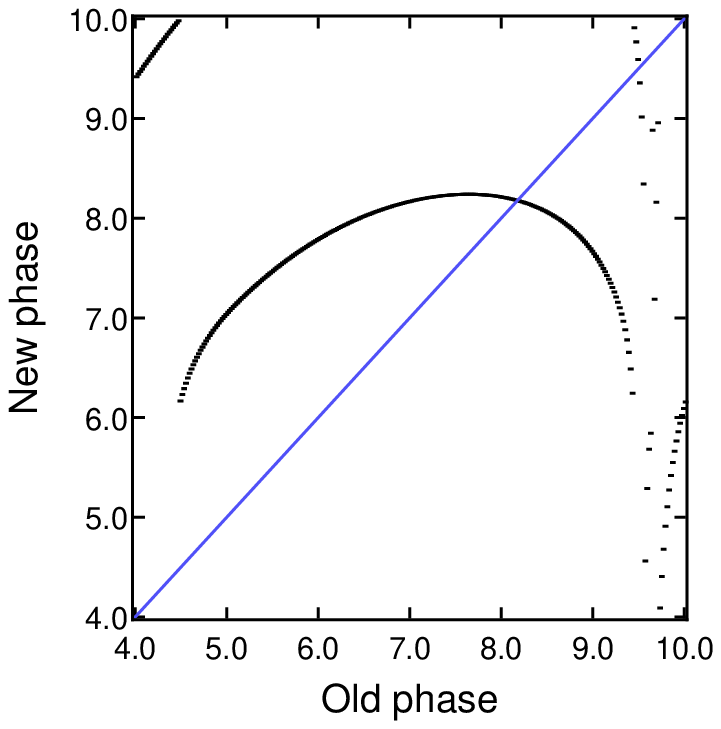}
\caption{The first return map $R_{f_T}$ to the interval $[4,10]$ (chosen
to enclose the plateau), for $A=10$ and $T=17.6$.  The blue line marks
the diagonal.}
\label{fig:first-return-map}
\end{figure}


\paragraph{Periodicity in drive period $T$.}
The approximate periodicity of $\lyap(F_T)$ seen in Figure
{\ref{fig:lyaps}} is easy to understand heuristically: kicking the
oscillator every $T$ seconds and kicking it every $T+T_0$ seconds should
yield the same asymptotic response because the oscillator simply
traverses $\gamma$ at frequency $1/T_0$ between kicks.  One can restate
this using phase resetting curves: if the drive period $T$ is
sufficiently large and $\theta(x_0)=t_0$, then the $f_T$-orbit $\of{t_0,
f_T(t_0), f_T^2(t_0), ...}$ should closely follow the phases
$\of{\theta(x_0),\theta(F_T(x_0)),\theta(F_T^2(x_0)), ...}$ of the
corresponding $F_T$-orbit.  Since $f_{T+T_0} = f_T$, this suggests that
$\lyap(F_{T+T_0})\approx\lyap(F_T)$.


\paragraph{Preference for entrainment.}
Figure {\ref{fig:prc1}} shows phase resetting curves for the {\HH}
equations for various values of drive period $T$ and drive amplitude
$A$.  For sufficientlly small values of $A$, the phase resetting curves
are circle diffeomorphisms: either there are sinks ({\em i.e.} stable
fixed pionts or stable periodic orbits), or the map is conjugate to a
rotation on a circle and the response of the kicked oscillator drifts
relative to the periodic drive.  As $A$ increases, the graph of $f_T$
rather quickly folds over and acquires critical points.  A striking
feature of the graphs in Figure {\ref{fig:prc1}} is the ``plateau,'' a
phase interval over which $f_T$ varies very slowly.  Another striking
feature is the ``kink'' around $\theta\approx 9.8$.  These features are
discussed in more depth in {\S\ref{sec:more-results}}.  For now, notice
that the plateau provides a simple mechanism for creating sinks:
changing the kick period $T$ shifts the graph of $f_T$ vertically.
Whenever the graph intersects the diagonal with a derivative $\abs{f'_T}
< 1$, then a stable fixed point is created.

\begin{table}
\begin{center}
\begin{tabular}{|c|c|c|}
\hline
Drive amplitude $A$ & Prob. of sink near plateau &
Prob. of $\lyap < 0$\\\hline
5 & 41\% & 48\%\\
10 & 58\% & 62\%\\
20 & 68\% & 70\%\\
30 & 76\% & 78\%\\\hline
\end{tabular}
\end{center}
\caption{Estimates of the probability of obtaining sinks near the
  plateau, as a function of $A$.  The data for this table is computed by
  trying about 40 values of $T$ for each choice of $A$ and examining the
  graph of the first return map to the interval $[4,10]$ (chosen to
  coincide with the ``plateau'') and its intersection(s) with the
  diagonal.}
\label{tab:prob-sinks}
\end{table}

This mechanism can be used to verify the results of Figure
{\ref{fig:lyaps-probs}}: compute the graph of the first return map of
$f_T$ to an interval around the plateau, then shift the graph vertically
using a number of different values of $T$ and estimate the fraction of
$T$'s for which $f_T$ has a stable fixed point (see Figure
{\ref{fig:first-return-map}}).  Table {\ref{tab:prob-sinks}} shows the
results.  For $A=10$, the $58\%$ probability of sinks corresponds fairly
closely with Figure {\ref{fig:lyaps-probs}}.  It is unclear whether the
ambiguous exponents in Figure {\ref{fig:lyaps-probs}} really represent
positive or negative Lyapunov exponents.  If a significant fraction of
the ambiguous exponents are really negative, then they must come from
small sinks.

\paragraph{Note on numerics.}
Phase resetting curves are computed here using a variation of the
Ermentrout-Kopell adjoint method {\cite{brown,ek}}.  The method is
described in the Appendix.  A systematic comparison of this method to
existing methods for computing phase resetting curves is beyond the
scope of the present paper and will be presented elsewhere.

\subsection{Response to repeated pulses: Wang-Young theory}
\label{sec:wy}

Phase resetting curves provide simple, intuitive explanations for many
dynamical properties of pulse-driven nonlinear oscillators.  For our
spiking {\HH} oscillator, explicitly-computed phase resetting curves
show why our pulse-driven neuron prefers entrainment over chaos.  In
order to infer asymptotic behavior, there needs to be a correspondence
between the orbits of $f_T$ and $F_T$, and the phase of
$x\in\Basin{\gamma}$ generally does not determine the phase of $F_T(x)$:
it may only do so approximately for a finite number of iterates.
When $\lyap(f_T) < 0$, this is enough to show that $f_T$ orbits indeed
approximate the phases of $F_T$.  Inferring chaotic behavior for $F_T$
from $f_T$ is far more difficult.  Wang-Young theory provides a
mathematical framework for inferring chaotic behavior using phase
resetting curves, and in addition explains why chaotic phenomena (and
all the other scenarios) is prevalent.

\begin{figure}
\begin{center}
\subgrapha{5in}{5in}{0.4}{}{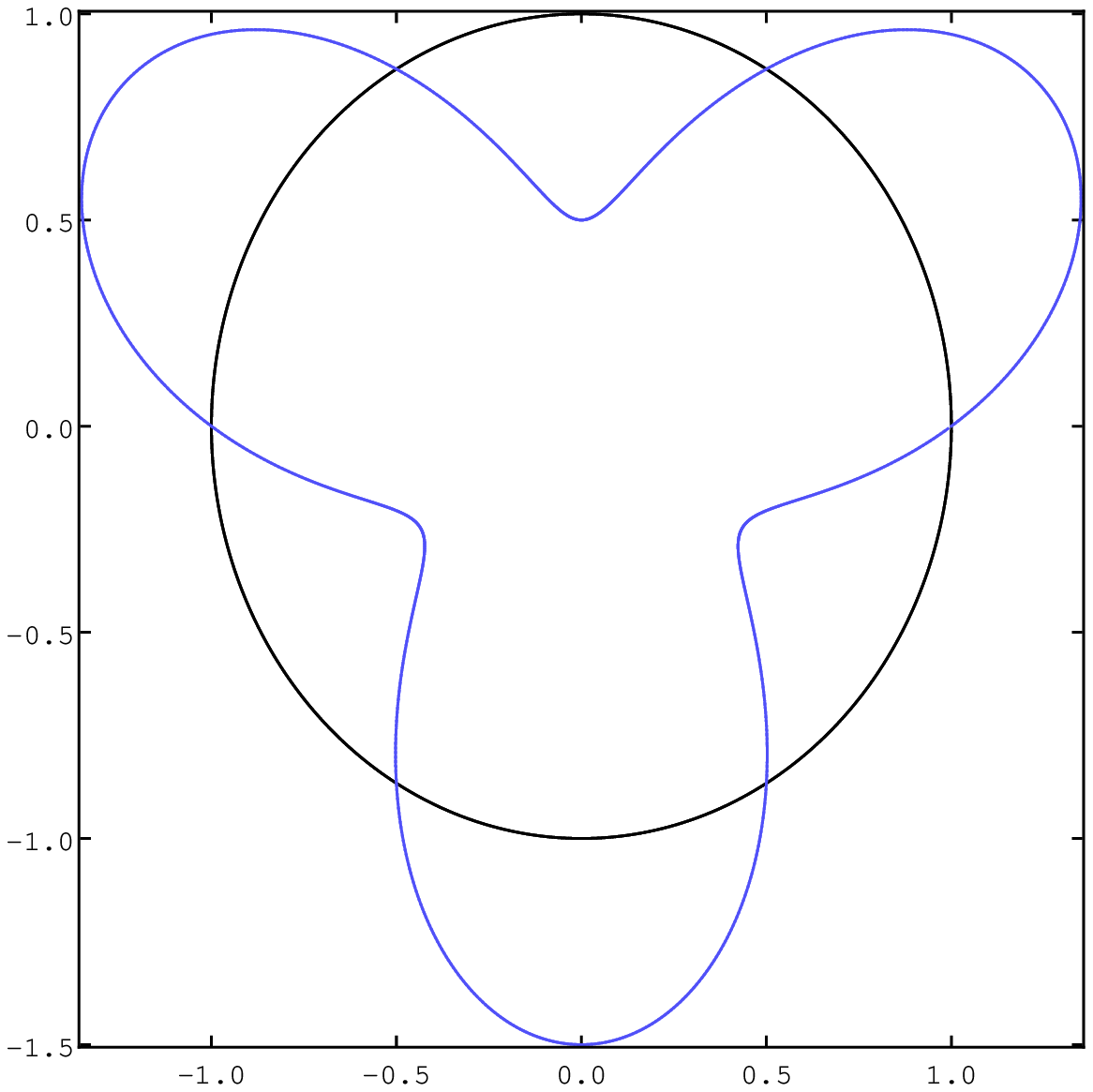}\qquad
\subgrapha{5in}{5in}{0.4}{}{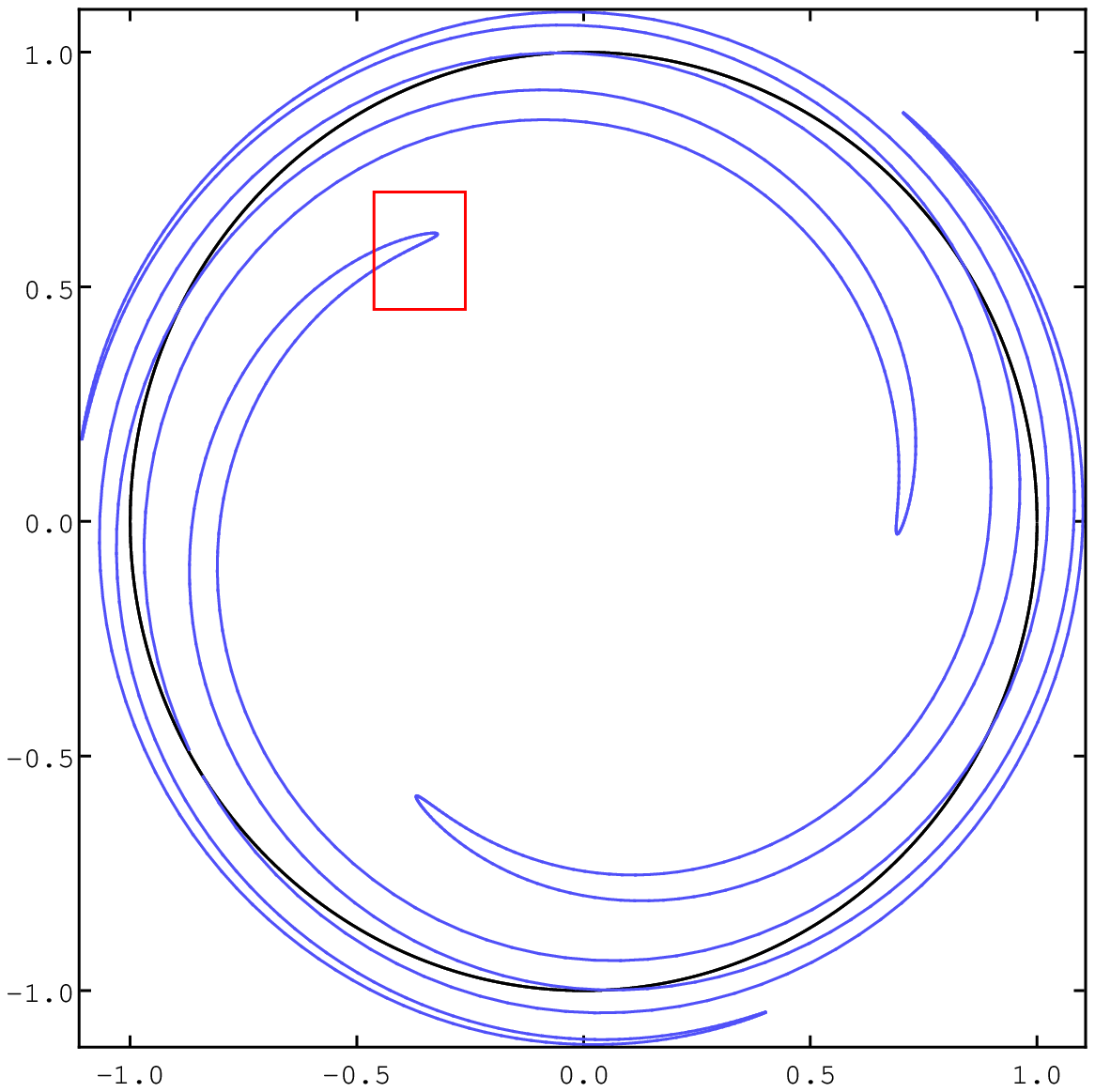}
\end{center}
\caption{A simple model (see Eq.~\ref{eqn:wy2}) which illustrates the
effect of shear.  In (a), a large number of initial conditions are
placed around a limit cycle (the black circle) and a radial kick is
applied to each point.  The blue curve shows the resulting positions of
each initial condition.  In (b), the kicked points are allowed to flow.
The red box contains one of the turning points (see {\S\ref{sec:wy}}).}
\label{fig:simple}
\end{figure}

{\em Shear} is an important ingredient of Wang-Young theory.  Let
$\gamma$ be the limit cycle which represents the unforced nonlinear
oscillator.  Near $\gamma$, the dynamics follows the periodic,
rotational motion on $\gamma$.  Shear refers to the presence of an
angular velocity {\em gradient} around $\gamma$: the stronger the shear,
the sharper the angular velocity changes at $\gamma$.  In two
dimensions, this means the flow runs much faster on one side of $\gamma$
than on the other; in the presence of strong shear, strong stable
manifolds tend to become more nearly tangent to $\gamma$.

Shear and its interaction with kicks is illustrated in a simple model in
Figure {\ref{fig:simple}}.  In the presence of strong shear, most ways
of kicking the oscillator which take advantage of shear ({\em e.g.}
kicks which do {\em not} take $x\in\gamma$ too close to the
strong-stable manifold $\varsub{W}{ss}(x)$) will cause segments of the
limit cycle to stretch and fold as they fall back toward $\gamma$.  The
phase space stretching caused by shear manifests itself in phase
intervals over which $\abs{f'_T} > 1$.  This expansion is conducive to
chaotic behavior.  However, shear also creates region of contraction
around ``turning points,'' an example of which is highlighted in Figure
{\ref{fig:simple}}.  Such turning points correspond to critical points
on the phase resetting curve and can easily counteract the expansion
needed for a positive Lyapunov exponent.  The competition between
expansion and contraction is the main source of difficulty in proving
$\lyap(F_T) > 0$.

To infer the existence of parameters for which the time-$T$ map $F_T$ is
fully chaotic, Wang and Young use results from their previous work on
strange attractors with $1$ expanding direction and a roughly toroidal
geometry {\cite{wy1,wy4}}.  Their main result gives conditions under
which there must be $T$ for which $\lyap(F_T) > 0$.  Furthermore, one
can find such ``chaotic parameters'' near ``nice'' values of $T$ for
which $f_T$ has a positive Lyapunov expoent.
Applying the general theory to kicked oscillators requires checking
certain geometric conditions.  This has been done for a few concrete
classes of models {\cite{wang,wy2,wy3}}.  To illustrate the consequences
of the theorems, consider the simple mechanical system {\cite{wy2}}
\begin{equation}
\label{eqn:wy2}
\ddot\theta(t) + \lambda\dot\theta(t) = \mu +
A\cdot\hat{K}(\theta(t))\sum_{n\in\Z}{\delta(t - nT)}.
\end{equation}
This model was first studied by Zaslavsky, who discovered that this
simple system can exhibit fully chaotic behavior {\cite{zaslavsky}}.
For Eq.~\ref{eqn:wy2}, it can be shown that the full range of scenarios
enumerated in the Introduction take place and that they are all
prevalent.  More precisely, Wang and Young prove (see Theorems 1-3 in
{\cite{wy2}}) that:
\begin{enumerate}

\item {\bf Invariant curve \& weak kicks:} When the drive amplitude $A$
  is sufficiently small (which is equivalent to having a large enough
  contraction rate $\lambda$), there exists a simple closed curve
  $\tilde\gamma$ to which all orbits of $F_T$ converge and which is
  invariant under $F_T$.  Moreover, we have the following dichotomy:
  \begin{enumerate}

  \item {\bf Quasiperiodic attractors:} There exist a set of $T$ of
  positive Lebesgue measure for which $F_T$ is topologically conjugate
  to an irrational rotation.  In this case, $F_T$ is uniquely ergodic on
  $\tilde\gamma$.

  \item {\bf Gradient-like dynamics:} There exists an open set of $T$
  such that $F_T$ has a finite number of periodic sinks and saddles on
  $\tilde\gamma$, and every orbit converges to one of these periodic
  orbits.

  \end{enumerate}

\item {\bf Gradient-like dynamics without an invariant curve:} 
   As $A$ increases (or $\lambda$ decreases), the invariant curve
  $\tilde\gamma$ breaks up.  Nevertheless, there continues to be an
  invariant set (no longer a simple closed curve) on which gradient-like
  dynamics persists.

\item {\bf Transient chaos:} For even larger $A$ or smaller $\lambda$,
  Smale horseshoes (see {\cite{gh}}) will form.  Horseshoes can coexist
  with sinks and saddles, creating transient chaos.

\item {\bf Chaos:} In the presence of sufficiently strong shear, there
  exists a positive measure set of drive periods $T$ for which $F_T$ is
  fully chaotic in the sense that it possesses (i) a strange attractor
  with a positive Lyapunov exponent; (ii) at least $1$ and at most
  finitely many ergodic SRB measures\footnote{SRB measures are natural
  invariant measures for dissipative dynamical systems.  They
  characterize the asymptotic behavior of a Lebesgue-positive measure
  set of initial conditions and have a number of nice mathematical
  properties.  See Young {\cite{young}} for an introduction.}  with no
  zero Lyapunov exponents; (iii) a central limit theorem; (iv)
  exponential decay of correlations if a power $F_T^N$ is mixing for
  some SRB measure $\nu$.
\end{enumerate}
Note that this list of (fairly well-understood) scenarios may not be
exhaustive.  Other scenarios or combinations of scenarios are not
excluded by the theory.  Also, the kicks in Eq.~\ref{eqn:wy2} are purely
radial.  This is not strictly necessary; any kick map which takes
advantage of shear will do.  See {\cite{wy2}} for precise conditions and
proofs.

\section{Further results}
\label{sec:more-results}

\subsection{More on the {\HH} phase resetting curve}

\subsubsection*{Plateau}


The plateau in the phase resetting curve for our pulse-driven {\HH}
model (see Figure {\ref{fig:prc1}}) corresponds a segment $\bar\gamma$
of the limit cycle $\gamma$ which becomes nearly parallel to a
strong-stable manifold after receiving a kick.  This can be seen by
examining the factors which contribute to the derivative $f'_T$ and
which can potentially cause $f'_T$ to become small over a relatively
large phase interval.  This can be checked by writing $f_T$ as a
composition of other functions and differentiating.

Let $\gamma:\R\to\R^4$ denote the limit cycle trajectory.  If we choose
$\gamma(0)$ so that $\theta(\gamma(0))=0$, then $\theta(\gamma(t))=t$
for all $t\in[0,T_0)$, and
\begin{align*}
f_T &=\theta\circ F_T\circ\gamma\\
&=\theta\circ\phi_T\circ K_A\circ\gamma.
\end{align*}
Changing $T$ does not affect $f'_T$, so we can set $T=0$.  Let $f=f_0$.
Then $f =\theta\circ K_A\circ\gamma$ and the chain rule gives
\begin{equation}
f' =\of{D\theta\circ K_A\circ\gamma}\cdot\of{DK_A\circ\gamma}\cdot
\dot{\gamma}.
\end{equation}
But $K_A(v,m,n,h)=(v+A,m,n,h)$, so its Jacobian $DK_A$ is the identity
matrix, and for all $t\in[0,T_0)$,
\begin{align}
\label{eqn:prc-factors}
f'(t) &= D\theta(K_A(\gamma(t)))\cdot\dot{\gamma}(t)\nonumber\\
&=\abs{D\theta(K_A(\gamma(t)))}\cdot\abs{\dot\gamma(t)}\cdot\cos(\mbox{angle}(D\theta(K_A(\gamma(t))),\dot\gamma(t)))\\
&=\abs{D\theta(K_A(\gamma(t)))}\cdot\abs{\dot\gamma(t)}\cdot\sin(\measuredangle(t)),\nonumber
\end{align}
where $\measuredangle(t)$ is the angle between $\dot{\gamma}(t)$ and the
strong-stable manifold at $K_A(\gamma(t))$.  The last step uses the fact
that the phase function $\theta:\Basin{\gamma}\to[0,T_0)$ is constant on
strong-stable manifolds (see {\S\ref{sec:theory}}).  This implies that
the gradient $D\theta(x)$ is everywhere orthogonal to
$\varsub{W}{ss}(x^*)$, where $x^*$ is the unique point in $\gamma$
having the same phase as $x$.  The factors in Eq.~\ref{eqn:prc-factors}
thus have simple, geometric meaning: $\gamma(t)$ is a point on the limit
cycle $\gamma$ and $\abs{\dot\gamma(t)}$ is the speed of the limit cycle
at that point; $\measuredangle(t)$ is the angle between $\dot\gamma(t)$
and the strong-stable manifold at $K_A(\gamma(t))$; and $\abs{D\theta}$
measures the rate at which the phase is changing at $K_A(\gamma(t))$.

\begin{figure}
\begin{center}
\subgraph{3in}{3in}{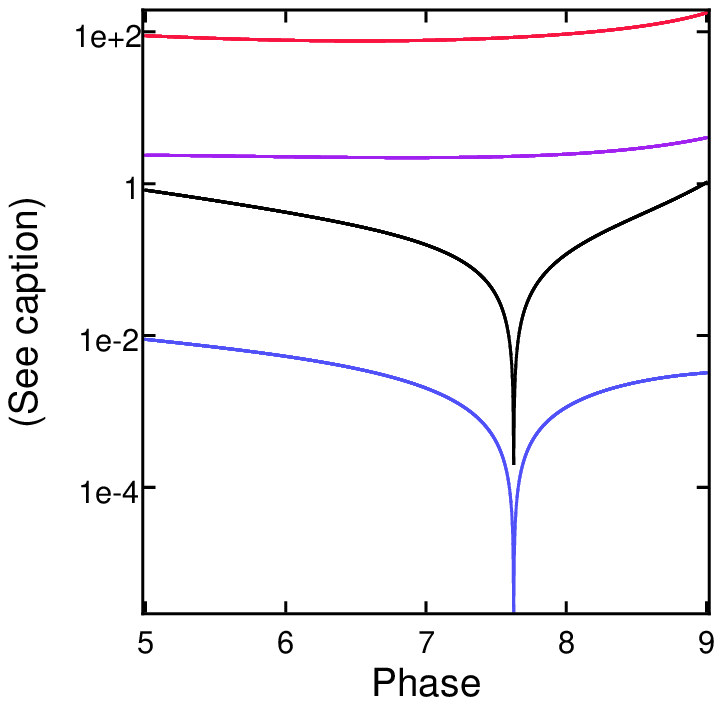}{}\qquad\subgraph{3in}{3in}{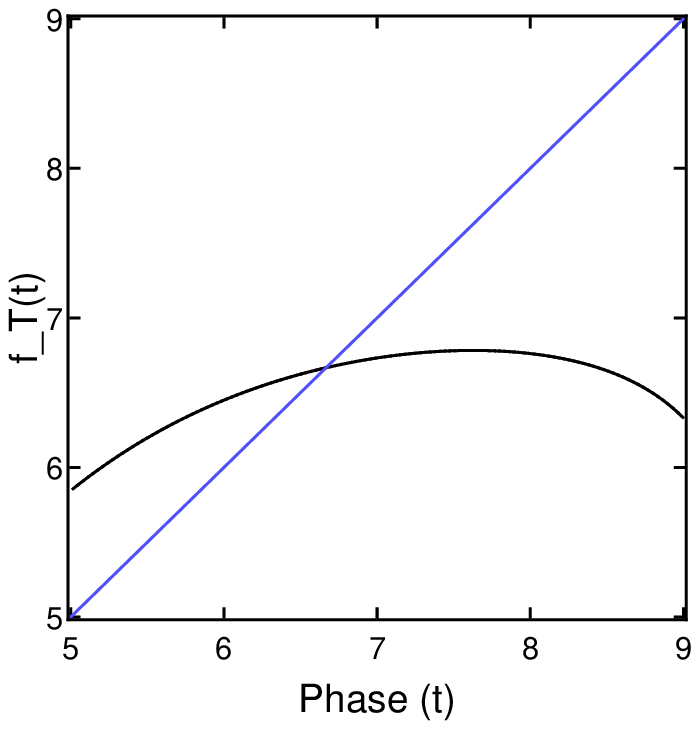}{}\\
\subgraph{4in}{3in}{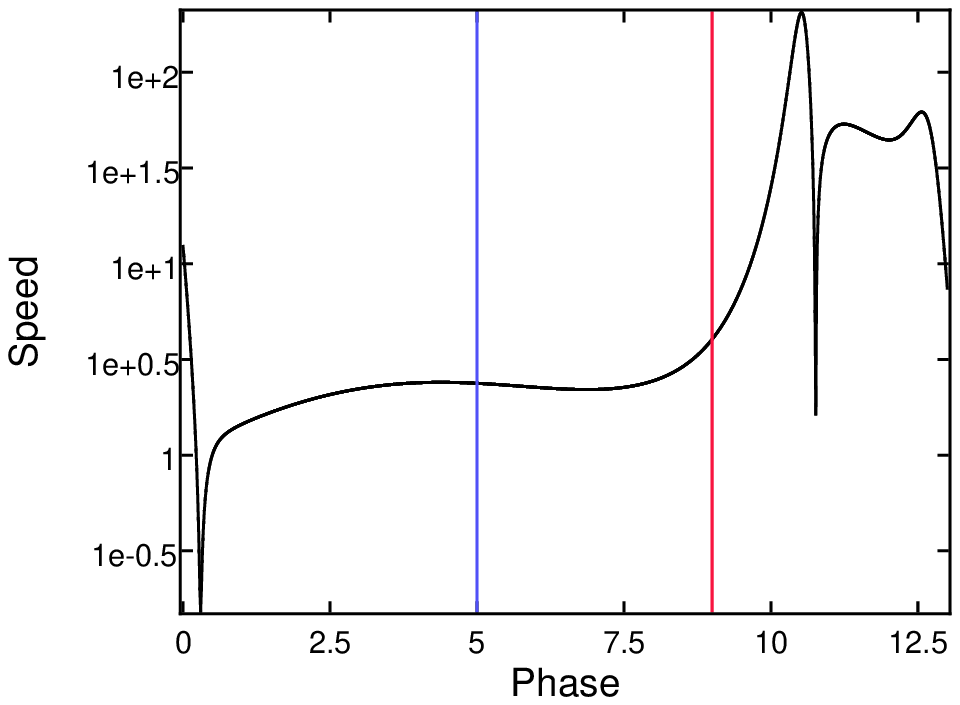}{}
\end{center}
\caption{The origin of the plateau: in (a), the curves are (i) {\em
Black:} $\abs{f'(t)}$; (ii) {\textcolor{blue}{\em Blue:}}
$\abs{\sin(\measuredangle(t))}$; (iii) {\textcolor{red}{\em Red:}}
$\abs{D\theta(K_A(t))}$; (iv) {\textcolor{purple}{\em Purple:}}
$\abs{\dot\gamma(t)}$.  In (b), the graph of $f_T$ near the plateau is
shown for reference.  Here, $A=10$.  (c) The speed of the {\HH} flow
along $\gamma$.  The vertical lines mark the interval $[5,9]$, which is
part of the plateau.}
\label{fig:plateau}
\end{figure}

Figure {\ref{fig:plateau}} shows $f'$ alongside the 3 factors in
Eq.~\ref{eqn:prc-factors}.  The figure shows that the plateau, where
$f'_T$ becomes nearly 0 over a long phase interval, coincides with the
near-vanishing of $\measuredangle(t)$.  The other factors of $f'$ stay
nearly constant over this interval.  Thus, there is a segment
$\bar\gamma$ of the limit cycle $\gamma$, corresponding to the phase
interval where $\measuredangle(t)$ is small, such that $K_A(\bar\gamma)$
is nearly tangent to a strong-stable manifold.  That the segment
$\bar\gamma$, which may be small as a subset of $\R^4$, corresponds to a
large phase interval, is due to the relatively slow speed of the limit
cycle near $\bar\gamma$.

What this argument does not explain is the robustness of this tangency
(equivalently, the robustness of the plateau) as the drive amplitude $A$
increases (see Fig. {\ref{fig:prc1}}).  This requires a detailed
analysis of the geometry (in $\R^4$!) of the strong-stable manifolds
(see {\S\ref{sec:future}}).

\paragraph{Numerical computation of $D\theta$.}
Figure {\ref{fig:plateau}} requires the numerical computation of the
gradient $D\theta(x)$ for $x\in\Basin{\gamma}$.  This can be done as
follows:

Fix $x\in\Basin{\gamma}$ and consider $\phi_{T_0}^n(x)$.  Clearly, the
limit $\lim_{n\to\infty}\phi_{T_0}^n(x)=x^*$ exists and has the property
that $x^*\in\gamma$, $\theta(x^*)=\theta(x)$, and
$x\in\varsub{W}{ss}(x^*)$.  Set $\projss(x) = x^*$.  Then $\projss$
projects $\Basin{\gamma}$ onto $\gamma$ and is the identity map on
$\gamma$.  Furthermore, the nullspace of the Jacobian matrix
$D\projss(x)$ of $\projss$ is the tangent to
$\varsub{W}{ss}(\projss(x))$ at $x$, by construction.

To compute $D\theta(x)$, the foregoing discussion suggests that we
compute $D\phi_{T_0}^n(x)$ for some large finite $n$.  For any finite
$n$, the singular values of $D\phi_{T_0}^n(x)$ consist of a dominant
singular value $\sigma_1$ and $3$ nearly zero singular values
$\sigma_2$, $\sigma_3$, and $\sigma_4$.  The $\sigma_i\to 0$ as
$n\to\infty$ for $i=2,3,4$.  Denote the left and right singular vectors
associated with $\sigma_1$ by $u$ and $v$.  It is easy to check that the
right eigenvector $v$ is orthogonal to the null space of $D\projss(x)$
and hence tangent to $D\theta$.

This computation requires a relatively accurate estimate of the
intrinsic period $T_0$ of the limit cycle $\gamma$, without which the
computation would not converge.  This paper adopts the following
strategy: instead of estimating $T_0$ just once and reusing its value,
solve the system of $24$ equations
\begin{equation}
\dot{x}_1 = H(x_1),\dot{x}_2 = H(x_2),\dot{J} = DH(x_2)\cdot J
\end{equation}
with initial conditions
\begin{equation}
x_1(0)\in\gamma,x_2(0)=x,J(0)=\mbox{Id}_{4\times 4}
\end{equation}
where $H$ is the {\HH} flow field and $x$ is the point at which we would
like to evaluate $D\theta$.  Note that $x_1, x_2\in\R^4$ and
$J\in\R^{4\times 4}$.  The solution of these equations then give $x_2(t)
=\phi_t(x)$ and $J(t) = D\phi_t(x)$.  The reference trajectory $x_1$ is
only used to count the number of periods which have elapsed, and the
trajectory $(x_2,J)$ is used to compute $D\projss(x)$.  This procedure
works fairly well in practice.

\subsubsection*{Kink}

\begin{figure}
\begin{center}
\subgraph{3in}{3in}{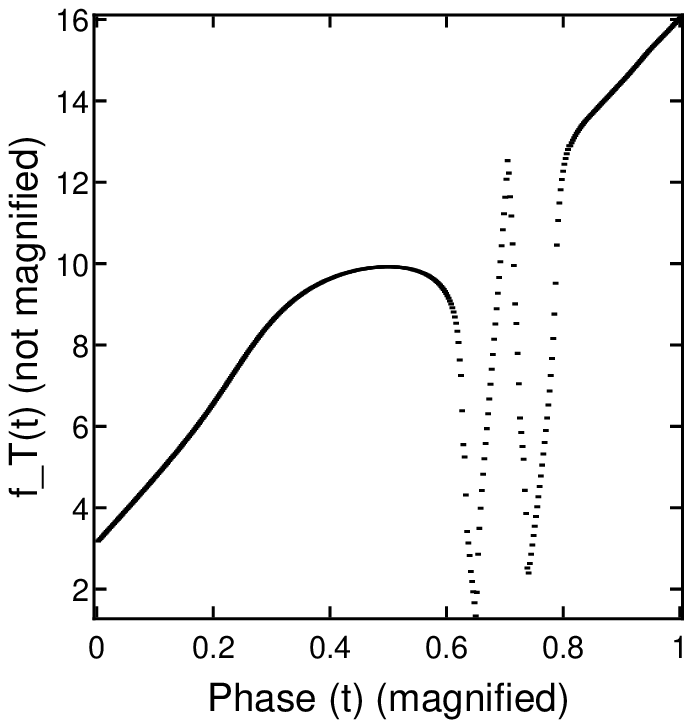}{}\qquad\subgraph{3in}{3in}{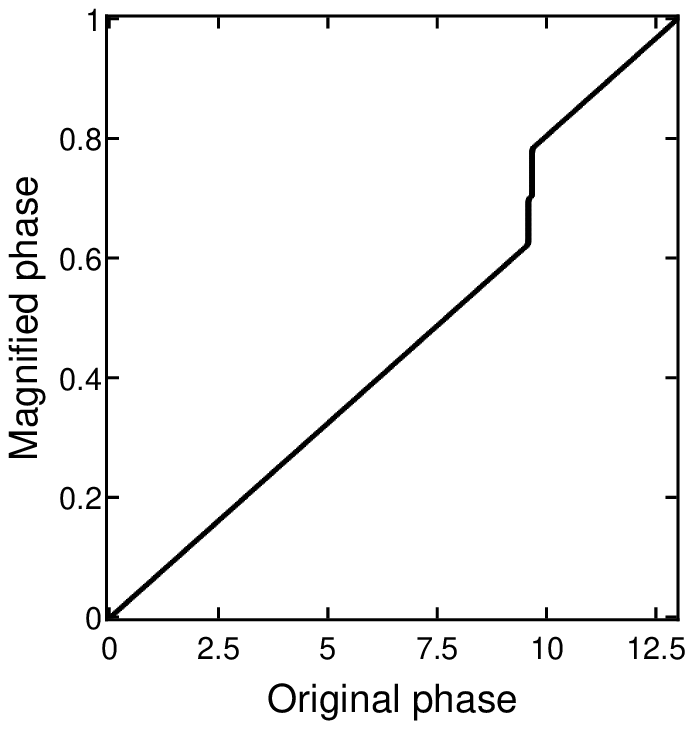}{}
\end{center}
\caption{(a) The graph of $f_T\circ g^{-1}$ with drive amplitude $A=10$,
with the abscissa shown in a new coordinate system $\theta' = g(\theta)$
to magnify the region around the ``kink.''  No interpolation is done in
this figure: only actually computed points are shown.  (b) The graph of
the coordinate transformation $g$.  The map $g$ is generated
automatically by the simple adaptive algorithm described in the
appendix.}
\label{fig:prc3a}
\end{figure}

\begin{figure}
\begin{center}
\subgraph{3in}{3in}{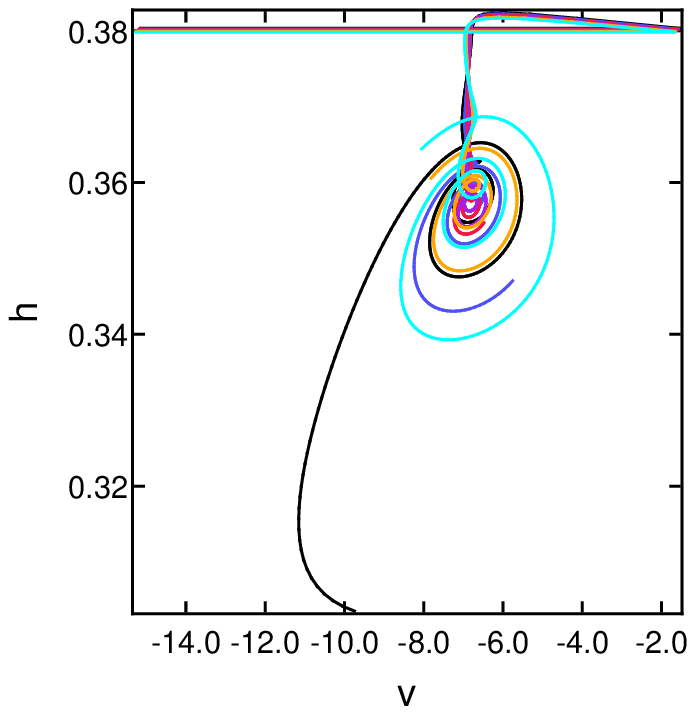}{}\qquad\subgraph{3in}{3in}{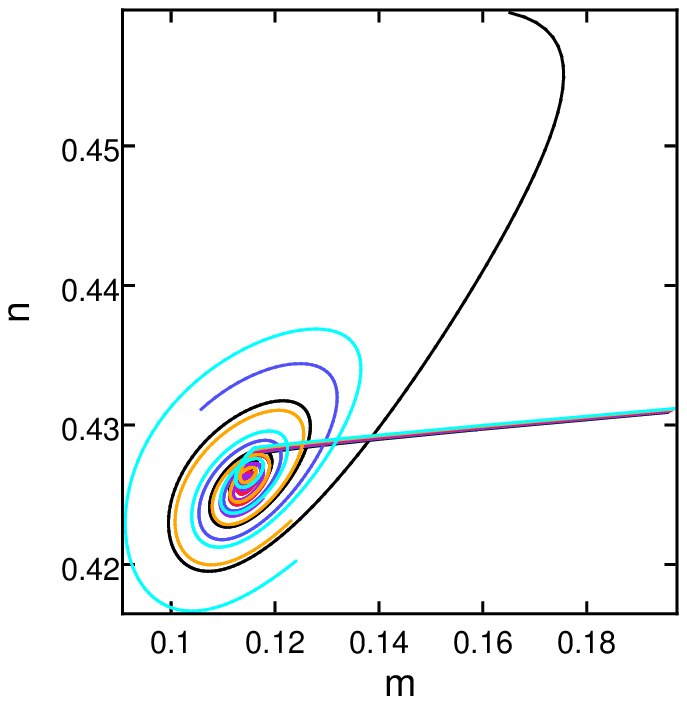}{}
\end{center}
\caption{A suggestive picture: in (a), a segment of $\gamma$, starting
in the upper left corner of the picture, is kicked straight across to
the upper right corner.  It then follows the flow toward the fixed point
for some time before spiraling away.  The overall direction of motion is
top to bottom.  (b) Another view of the approach to the fixed point.
The overall direction of motion here is right to left.  The kick
amplitude is $A=13.589$.}
\label{fig:kink}
\end{figure}

\begin{figure}
\begin{center}
\subgrapha{5in}{5in}{0.4}{}{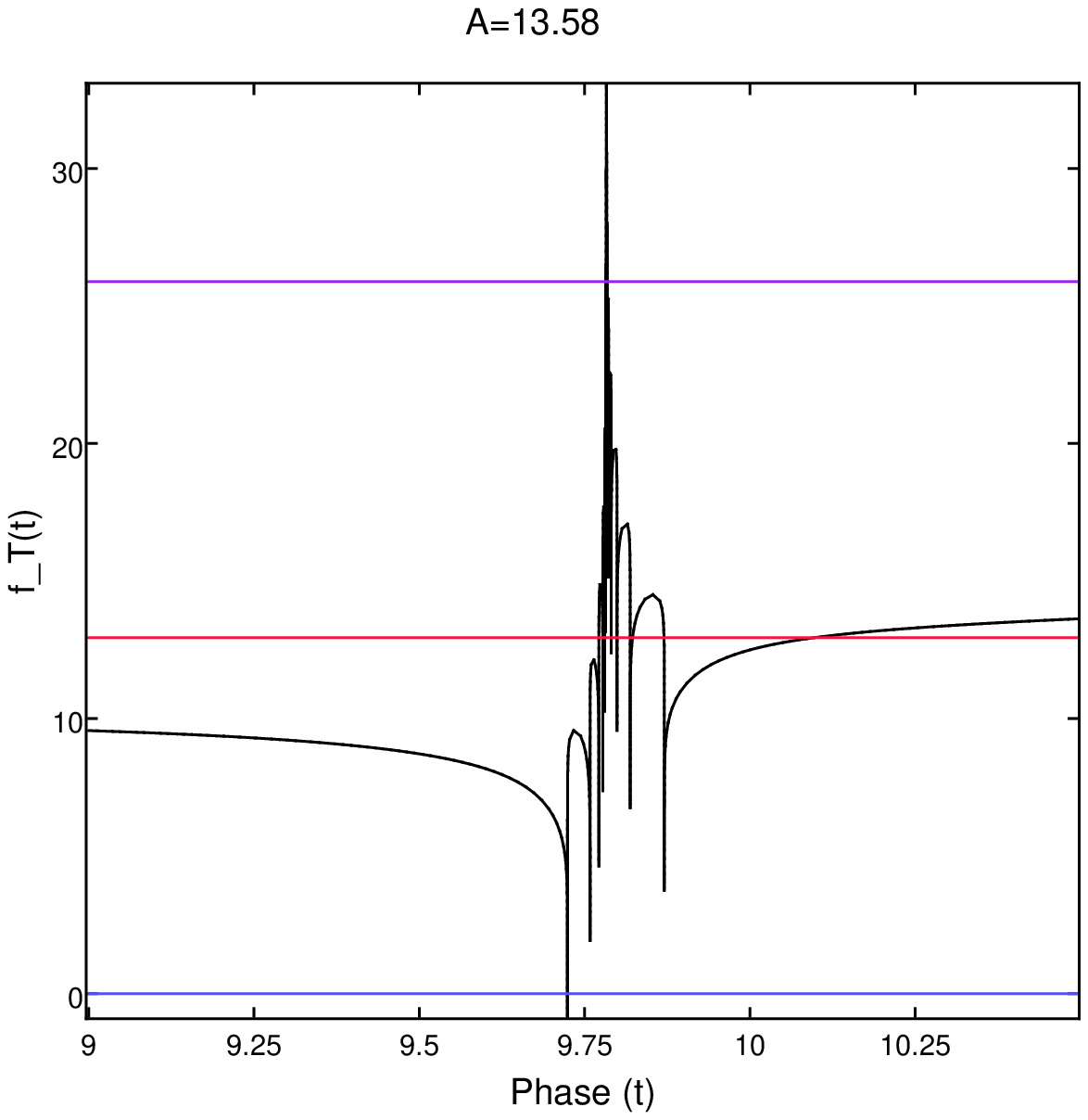}
\subgrapha{5in}{5in}{0.4}{}{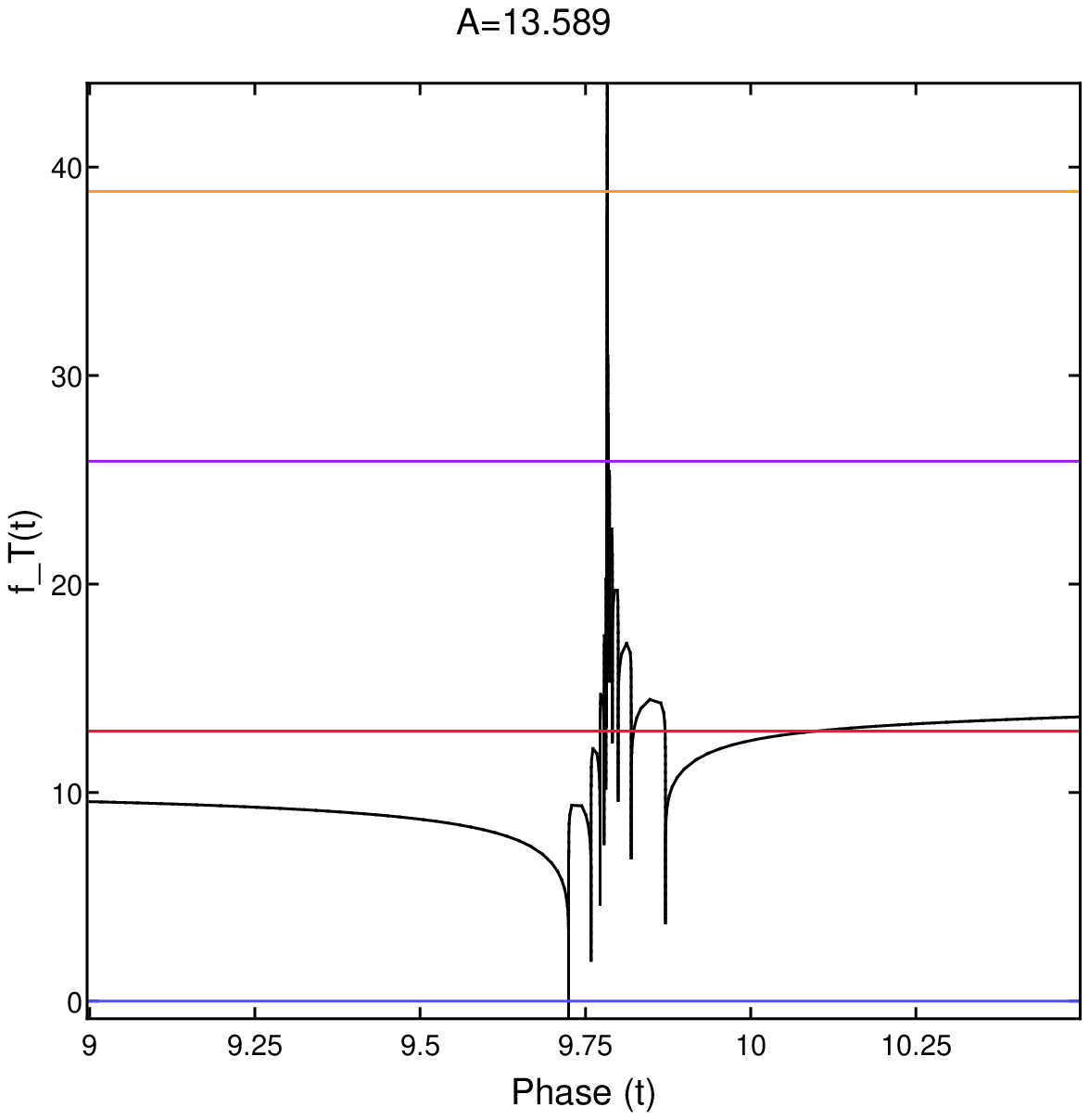}
\subgrapha{5in}{5in}{0.4}{}{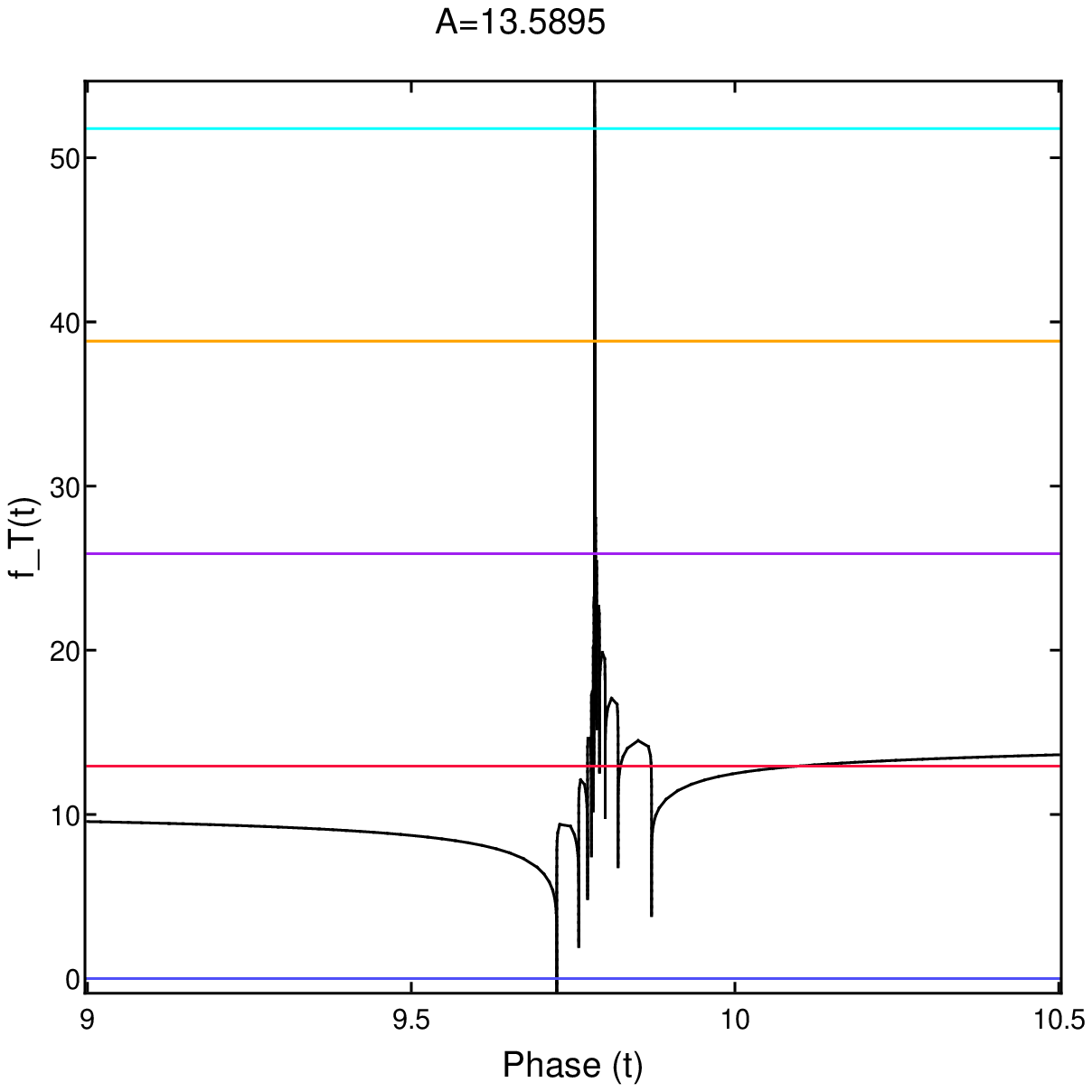}\\
\subgrapha{5in}{5in}{0.4}{}{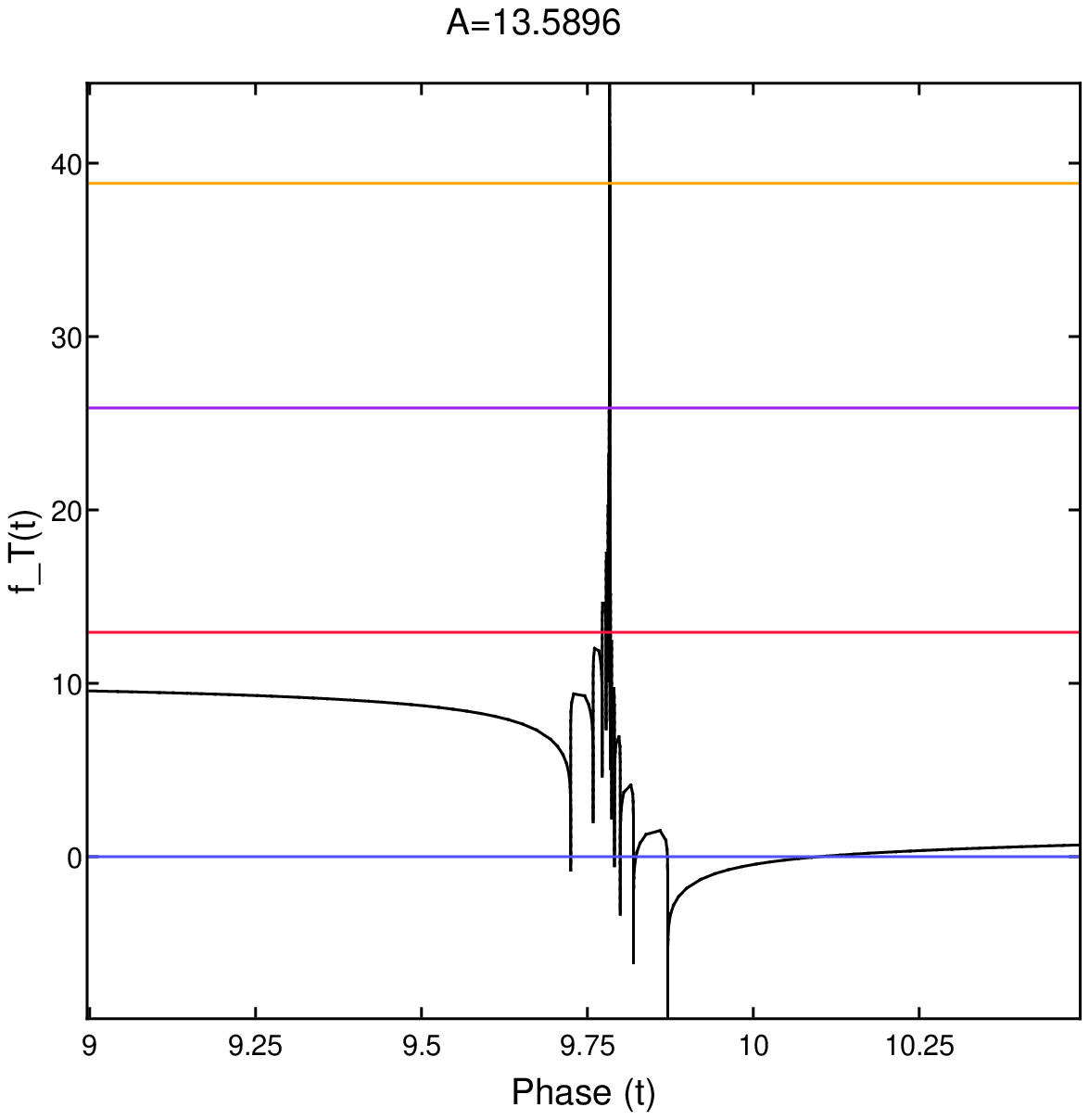}
\subgrapha{5in}{5in}{0.4}{}{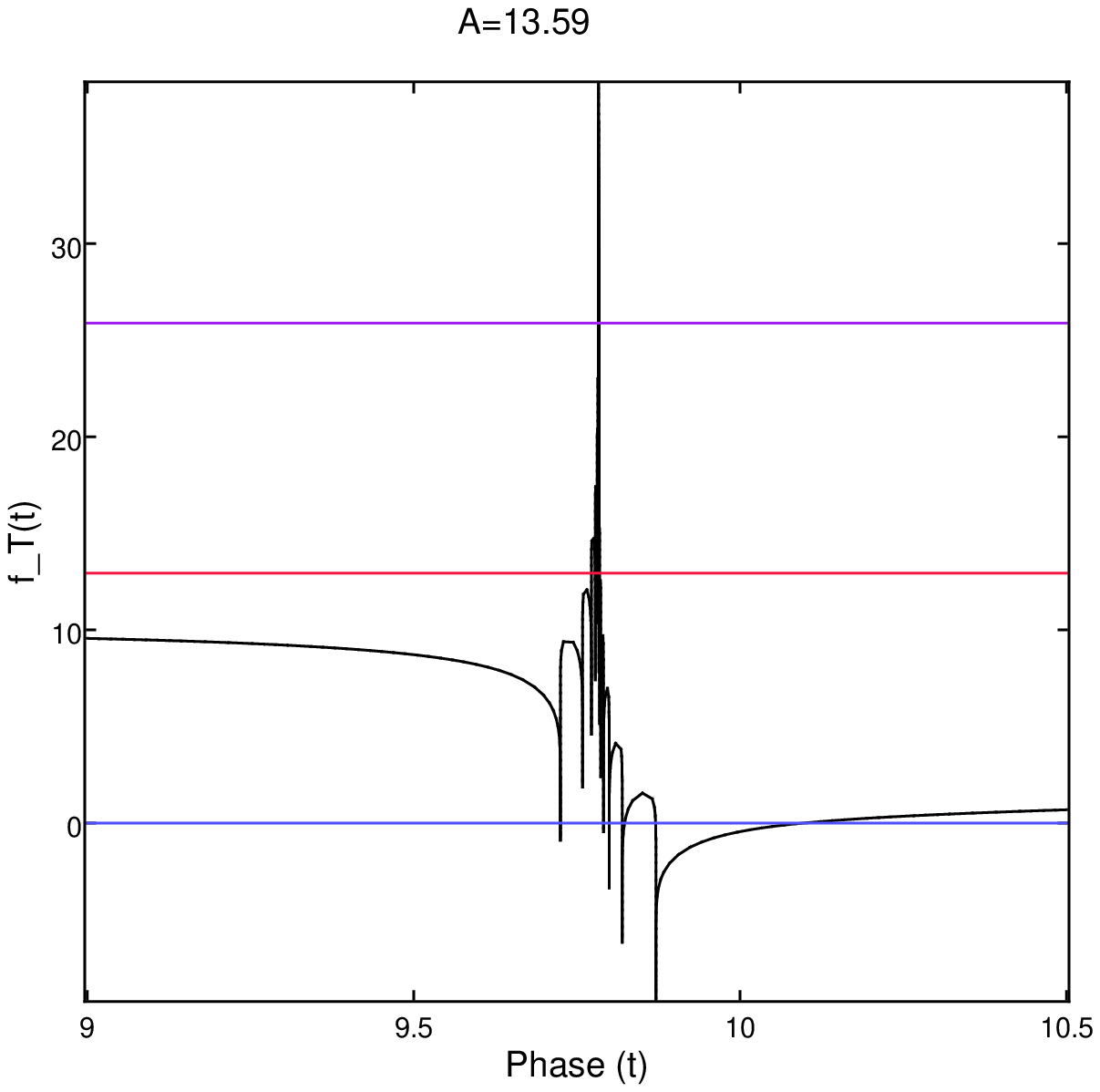}
\subgrapha{5in}{5in}{0.4}{}{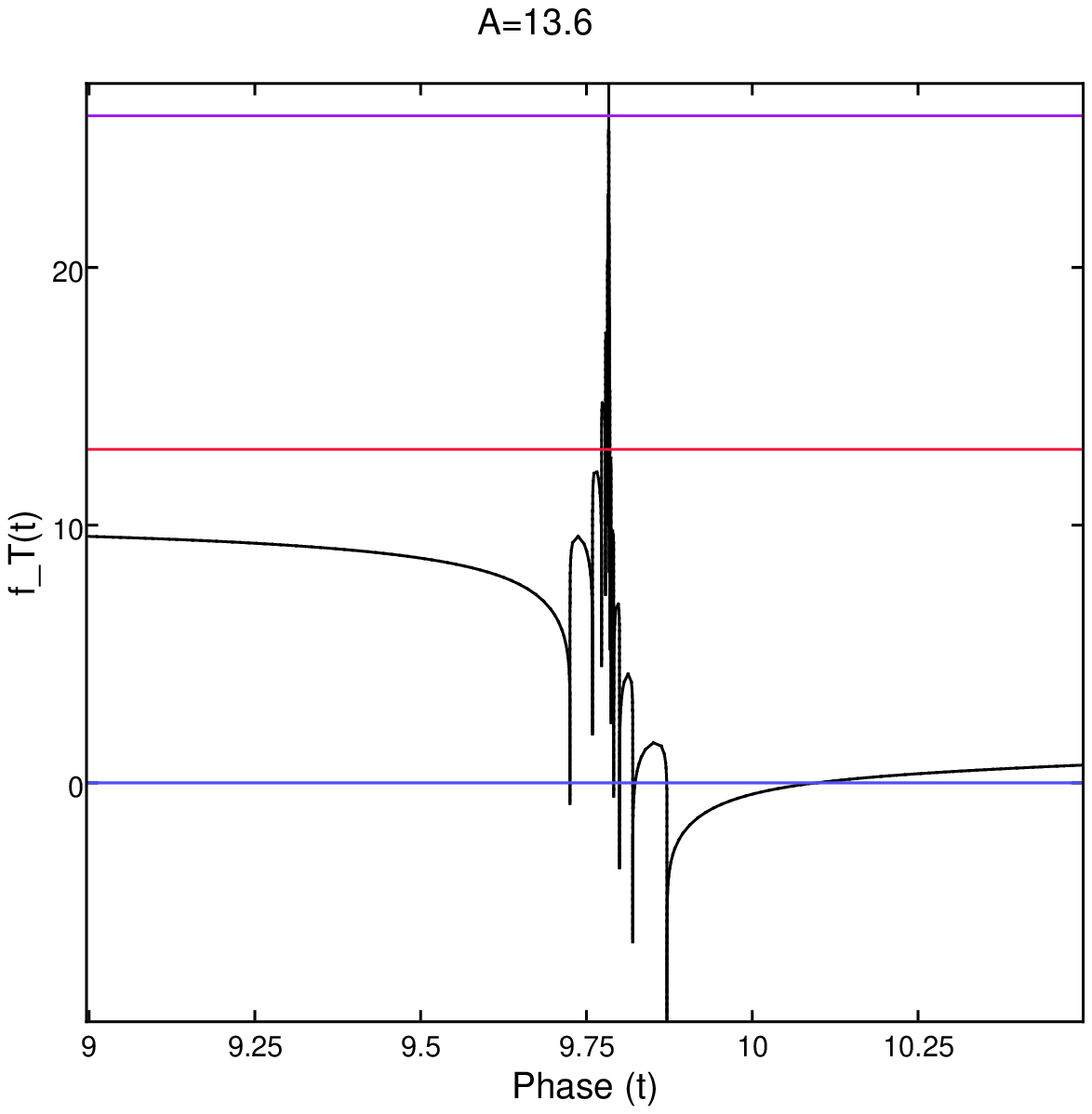}
\end{center}
\caption{Graphs of $f_T$ for kick amplitudes $A$ approaching
$\varsub{A}{crit}$ from below (a-c) and above (d-f).  The horizontal
lines mark integral multiples of $T_0$.}
\label{fig:kink1}
\end{figure}

It is natural to ask whether $f_T$ (see Fig. {\ref{fig:prc1}}), for
$A=10$ or $A=20$, is discontinuous around the kink.  A discontinuity
indicates that there are points in a neighborhood of $\gamma$ which can
be kicked outside the basin of $\gamma$.  This is not likely the case:
Figure {\ref{fig:prc3a}} shows a magnified view of the phase resetting
curve near the kink; the graph does not include any numerical
interpolation.  The figure is obtained by fixing a small parameter
$\delta > 0$ and adaptively refining the grid $\set{\theta_n}$ on which
the phase resetting curve is evaluated until
$\abs{f_T(\theta_{n+1})-f_T(\theta_n)}\leq\delta$.  In Figure
{\ref{fig:prc3a}}, $\delta$ is set to $0.1$.  The adaptive procedure
(see the Appendix) continues to converge for smaller values of $\delta$.

Figure {\ref{fig:kink}} suggests an explanation for the kink: that it is
likely caused by a segment of $\gamma$ being kicked near the stable
manifold of the unstable fixed point.  This would cause the segment to
wind around the stable manifold and eventually spirals away from the
fixed point.  (Recall that the two unstable eigenvalues of the fixed
point form a complex conjugate pair.)  In the process the kicked segment
spreads apart and its subsets pick up different amounts of time delays.
However, because the {\HH} phase space is $4$-dimensional, Figure
{\ref{fig:kink}} cannot give a reliable picture of the dynamics:
projecting onto $2$ dimensions loses too much information.

The scenario sketched above predicts that there exist a critical kick
amplitude $\varsub{A}{crit}$ at which $K_A(\gamma)$ intersects the
stable manifold of the fixed point.  (There may be more than 1
intersection, and more than 1 value of $A$ which cause intersections.)
As $A\to\varsub{A}{crit}$, the phase resetting curve should start
winding around $S^1$ more and more.  This can be numerically tested: an
estimate of $\varsub{A}{crit}$ is computed using the Nelder-Meade
algorithm {\cite{recipes}} to minimize the closest distance of a
trajectory to the fixed point.  This yields a critical value
$\varsub{A}{crit}\approx 13.58953...$.  When $A =\varsub{A}{crit}$
exactly, $f_T$ should wind around infinitely many times and possess a
singularity near the location(s) of intersection.  For
$A\neq\varsub{A}{crit}$, $f_T$ remains smooth, but as
$A\to\varsub{A}{crit}$, $f_T$ should develop a singularity and blow up.
See Figure~{\ref{fig:kink1}}.

\subsubsection*{Horseshoes \& transient chaos}
\label{sec:horseshoe}

\begin{figure}
\begin{center}
\putgraph{3in}{4in}{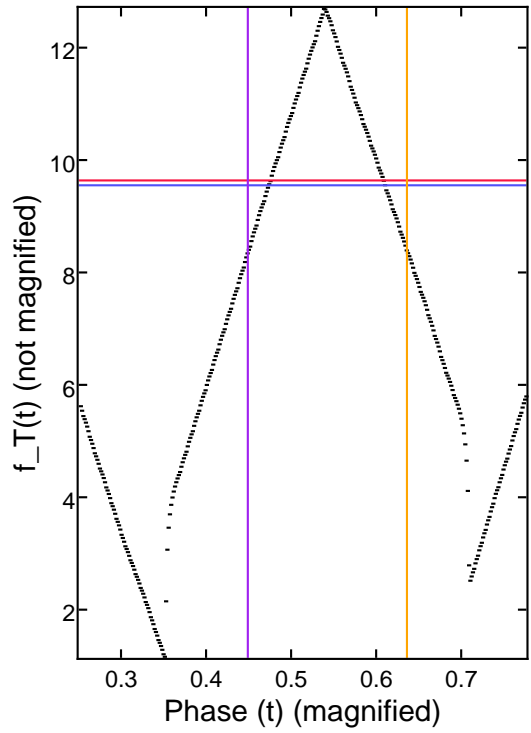}
\end{center}
\caption{The phase resetting map and an interval $I$ (marked by the
straight lines) which maps across itself twice.  The abscissa (but not
the ordinate) is shown in transformed coordinates.  Here, $A=10$ and
$T=81$.  This is a ``small'' horseshoe: because the derivative of $f_T$
is so large there (on the order of $10^3\sim 10^4$), most numerically
computed orbits escape the horseshoe after a few iterates.}
\label{fig:horseshoe}
\end{figure}

Wang-Young theory also guarantees the existence of $T$'s for which $F_T$
exhibits transient chaos, {\em i.e.} $F_T$ possesses a Smale horseshoe
{\cite{gh}} together with a sink.  The coexistence of a horseshoe with a
sink has the following effect on the dynamics: almost every $F_T$-orbit
would eventually fall into a sink, but an orbit which wanders near a
horseshoe would dance around unpredictably for a finite number of
iterations.  Two nearby orbits which enter the vicinity of a horseshoe
can emerge widely separated, and fall into the sink out of phase (unless
the sink happens to be a fixed point).  In terms of time series data,
this kind of behavior can be recognized by looking at {\em pairs} of
trajectories and finding that they chaotically ``flutter'' about before
settling down into a steady periodic motion, likely out of phase.


In contrast to entrainment and chaos, transient chaotic behavior is
difficult to observe in the pulse-driven {\HH} system.  This is because
the most likely place to find a horseshoe is near the kink, where the
expansion so strong that most trajectories escape very quickly.
Nevertheless, it is possible to find indirect evidence for horseshoes in
the pulse-driven {\HH} model.  To do so, one looks for an interval
$I\subset[0,T_0)$ such that $f_T(I)$ gets mapped completely across $I$
at least twice.  It is easy, for example, to find a ``small'' horseshoe
around the kink in the phase resetting curve.  See Figure
{\ref{fig:horseshoe}}.  The phase interval $I$ tells us the rough
location of a horseshoe for $F_T$.

To go from such an interval $I$ to a horseshoe for the full map $F_T$,
it is necessary to (i) blow up the corresponding segment of $\gamma$ to
form an open set $U\subset\R^4$ such that $F_T(U)$ intersects $U$ at
least twice, and the intersection stretches all the way across $U$ in
the unstable direction (along $\gamma$); and (ii) find invariant cones
{\cite{gh}}.  This can be done in a straightforward manner and is not
discussed further here.







\subsection{Miscellaney}
\label{sec:other-params}

\subsubsection*{Decay of correlations}

\begin{figure}
\begin{center}
\subgraph{3in}{3in}{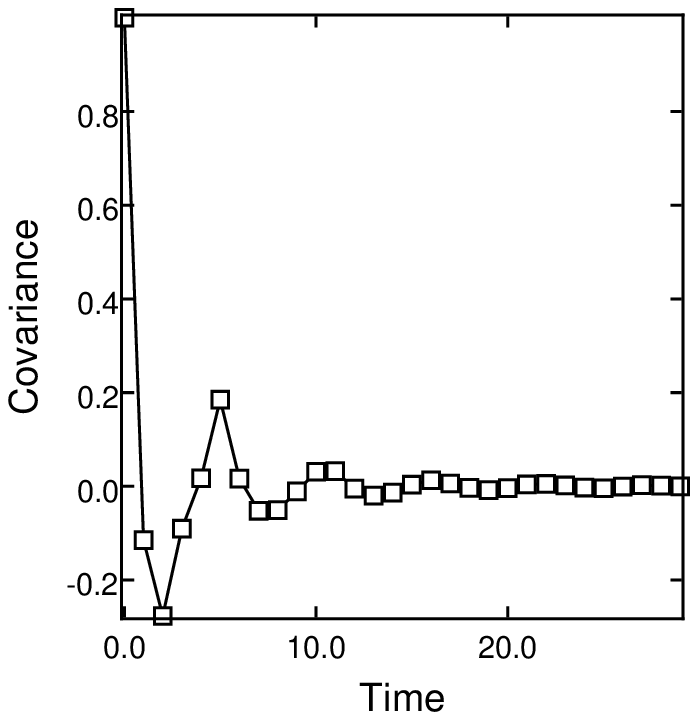}{Linear}\qquad
\subgraph{3in}{3in}{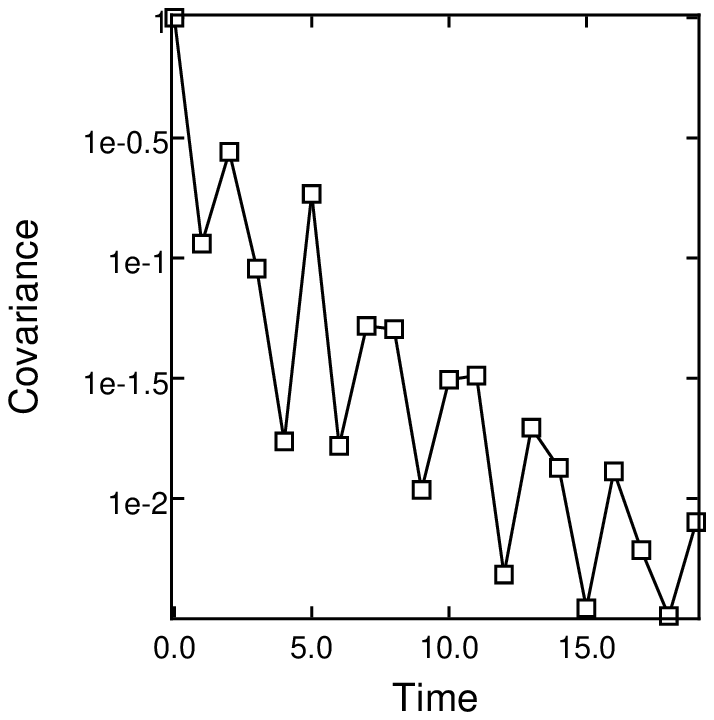}{Semilog}
\end{center}
\caption{Normalized autocorrelation function $C_{vv}(n)$ for the $v$
variable of the time-$T$ map $F_T$.  This graph indicates that the
dynamics is mixing, but the data is insufficient to confirm that the
system is exponentially mixing.}
\label{fig:corr}
\end{figure}

Wang-Young theory predicts that when the dynamics of a pulse-driven
oscillator is chaotic and there is a unique SRB measure, time
correlations (more precisely time autocovariance functions) decay
exponentially fast in time.  Figure {\ref{fig:corr}} shows the time
autocovariance function
\begin{equation}
C_{vv}(n) =\int{(v\circ F_T^n)\cdot v\ d\mu} -\of{\int{v\ d\mu}}^2
\end{equation}
for the voltage variable $v$ ($\mu$ is an ergodic invariant measure).
While it clearly decays as $n\to\infty$ and thus provides evidence that
the invariant measure is mixing, the data is not sufficient to confirm
that the decay is exponential.

\subsubsection*{Response to finite-duration pulses}

A natural variation on the numerical experiments of previous sections is
to replace instantaneous impulses with finite-duration pulses.
Heuristically, if the pulse durection $t_0$ is less than the fastest of
the intrinsic timescales of $\gamma,$ the resulting response should be
essentially the same as the response to instantaneous impulses.  With
$I\approx\Imain$, these time scales are $12.944$ (= the period), $5.1$,
$0.50$, and $0.12$ (corresponding to the negative Lyapunov exponents).

\begin{figure}
\begin{center}
\subgraph{3in}{3in}{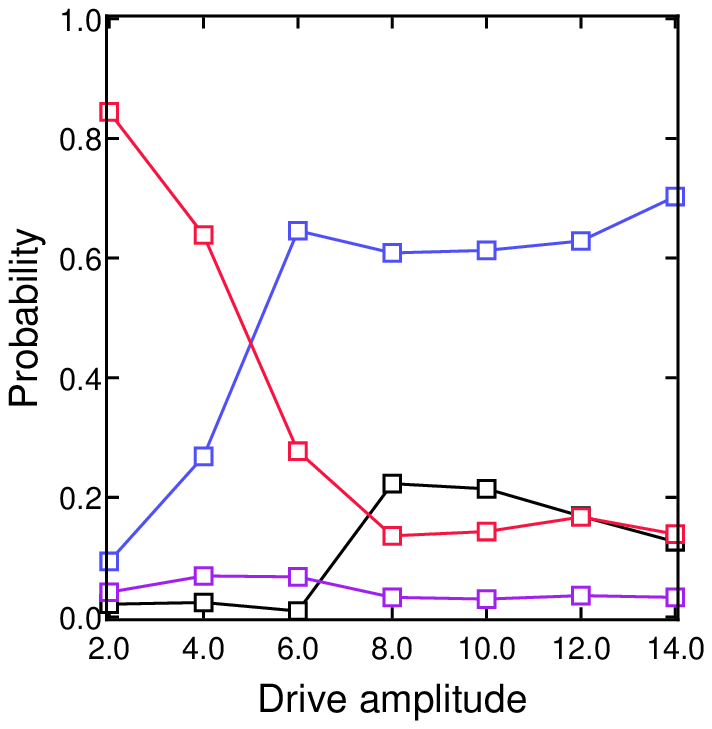}{}\qquad\subgraph{3in}{3in}{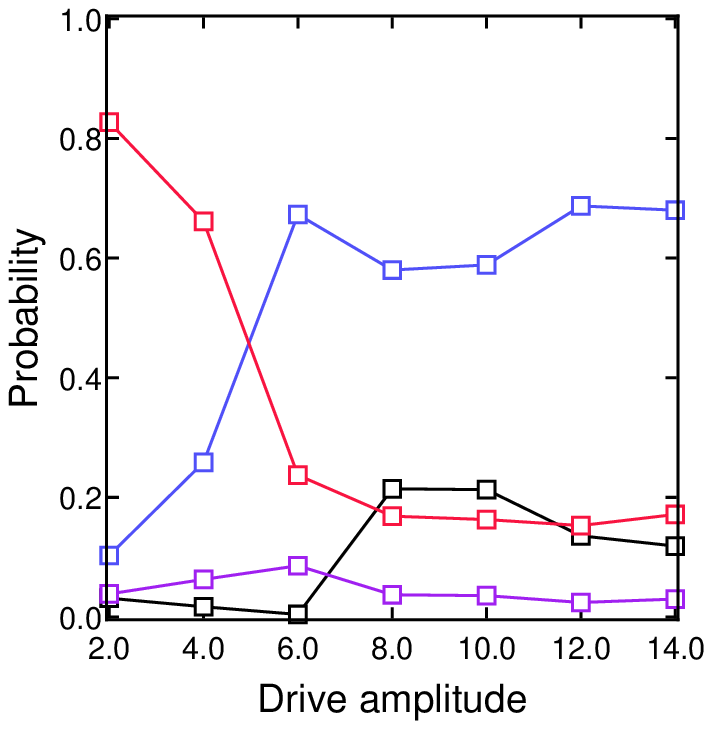}{}\\
\subgraph{3in}{3in}{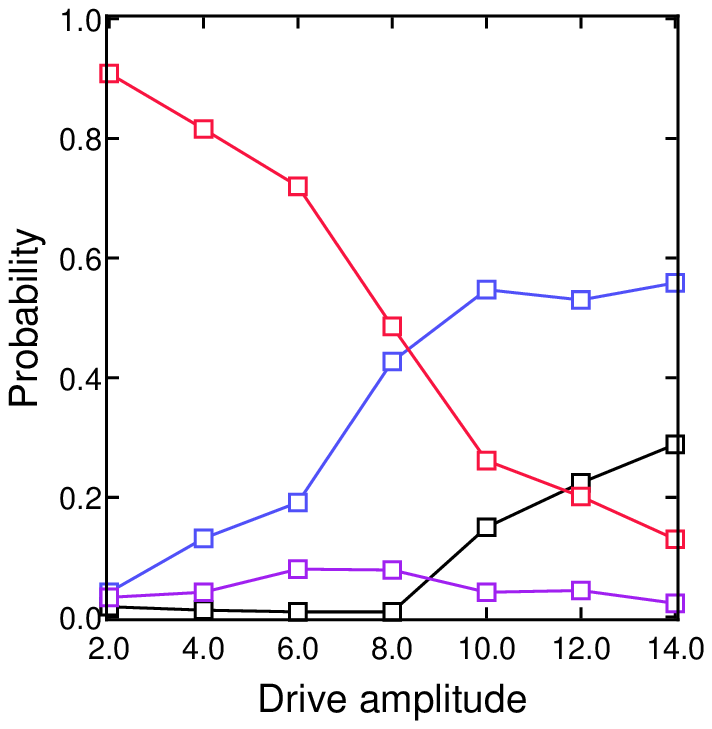}{}\qquad\subgraph{3in}{3in}{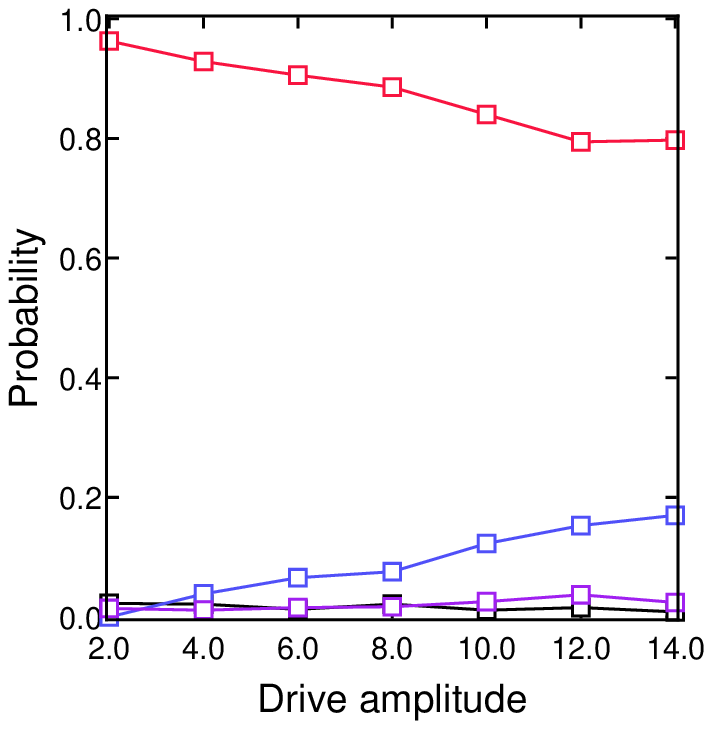}{}\\
\end{center}
\caption{Distribution of Lyapunov exponents for $I\approx\Imain$ with
  finite-duration pulse of duration (a) $t_0=0.05$, (b) $t_0=0.3$, (c)
  $t_0=2.75$, and (d) $t_0=9.0$.  See Figure {\ref{fig:lyaps-probs}}
  caption for details.  Recall that the Lyapunov exponents of the
  unperturbed limit cycle $\gamma$ are $-0.196$, $-2.01$, and $-8.31$,
  corresponding to relaxation times of $5.1$, $0.50$, and $0.12$.}
\label{fig:pulse}
\end{figure}

Figure {\ref{fig:pulse}} summarizes the numerical results for $t_0 =
0.05$ (shorter than all time scales), $0.3$ (shorter than all but one
time scale), $2.75$ (shorter than all but two fastest time scales), and
$9.0$ (very slow, not really pulsatile in any sense of the word).  These
graphs should be compared to Figure {\ref{fig:lyaps-probs}}.  The pulse
amplitude is adjusted so that the {\em total} amount of charge delivered
is the same as an impulse of amplitude $A$.  Interestingly enough, the
behavior seen earlier are quite robust and disappear only when $t_0 =
9$.  These results suggest that the contracting directions do not mix
very much over $\gamma$, and only the slowest contracting time scale
participates in the production of chaotic behavior.

\section{Outlook}
\label{sec:future}

The results reported here show that the pulse-driven {\HH} model
(\ref{eqn:pulse-driven-hh}) responds to low-frequency (relative to the
intrinsic period $T_0$ of the spiking neuron) periodic impulse forcing
in a wide range of ways.  Depending on the drive period and drive
amplitude, the response can range from entrainment to fully chaotic
behavior.  This is consistent with the predictions of Wang-Young theory.
Furthermore, as shown in {\S\ref{sec:theory}}, it is possible to explain
some phenomena specific to our pulse-driven {\HH} oscillator in terms of
special features of the phase resetting curve and provide a partial
understanding of the source of these features.

There are some interesting directions for future work:

\paragraph{Random kick times.}
The shape of the {\HH} phase resetting curve suggests that if one were
to drive a {\HH} neuron using a pulse train with {\em random} kick
times, the resulting random dynamical system can have a negative
Lyapunov exponent.  This is because the phase resetting curve moves up
and down from kick to kick, and for any kick time distribution which is
sufficiently uniform ({\em e.g.} an exponential distribution), the
probability that the plateau intersects the diagonal is high.  The size
of the plateau suggests that over many iterates, contraction may
dominate expansion, leading to a negative Lyapunov exponent.  A negative
Lyapunov exponent implies that two {\HH} neurons, when driven by a {\em
common} pulse train with random kick times, will synchronize.  That is,
the plateau provides a way to create a ``random fixed point''
{\cite{ledrappier-young}}.  These predictions are consistent with
preliminary numerical results and with a perturbation theory developed
by Nakao {\em et.~al.}  for randomly-kicked oscillators in the limit of
weak kicks {\cite{nakao}}.  The heuristic, geometric argument sketched
above may lead to an extension of their result to the regime of strong
kicks.

This synchronization mechanism has also been studied numerically by Doi
in the context of a simple pieceise linear map {\cite{doi1}}.  In
addition, there is an extensive literature on noise-induced synchrony in
neural models, including white-noise-driven {\HH} equations
{\cite{pakdaman,zhou-kurths}}.  Models exhibiting noise-induced
synchrony provide a concrete framework for exploring neural reliability
{\cite{hunter-milton,mainen}}.

\paragraph{Robustness of the phase resetting curve.}
How robust are the features (plateau, kink) of the {\HH} phase resetting
curve under perturbations to parameters?  How robust is the geometry of
the near-tangency of kicked segments and strong-stable manifolds
responsible for forming the plateau?  As the range of phenomena
predicted by Wang-Young theory may be present in two-dimensional models
like the Morris-Lecar or FitzHugh-Nagumo, these models may provide a
good starting point for exploring these questions in lower-dimensional
settings.

\appendix

\section{Appendix: Computation of phase resetting curves}
\label{sec:prc}

All numerical calculations of phase resetting curves reported in this
paper use the algorithm described here.  It is closely related to an
algorithm due to Ermentrout and Kopell {\cite{ek}}.  It is presented
here for the sake of completeness; a systematic comparison with existing
methods will appear elsewhere.  The algorithm computes the phase
resetting curve for {\em finite-size perturbations} but can be adapted
to compute phase resetting curves for infinitesimal perturbations.  See
Appendix A in {\cite{brown}} and
{\cite{oprisan-canavier,williams-bowtell}}) for more general discussions
of phase resetting curves, and Ermentrout, Pascal, and Gutkin
{\cite{ermentrout1}} for a discussion of computing phase resetting
curves experimentally.

The basic idea is to numerically compute the strong-stable linear
subspaces along the limit cycle.  Then, using these linear subspaces as
approximations to strong-stable submanifolds, project a kicked point
down to the limit cycle, where the phase can be estimated.  The
algorithm really computes strong-stable linear subspaces along $\gamma$
and uses these linear subspaces to approximate the phase resetting
curve.


Some preliminaries: the limiting map $\bar{F}_T
=\lim_{n\to\infty}F_{T+nT_0}$ can be characterized abstractly by the
equation
\begin{equation}
\label{eqn:phase-resetting-def-too}
\bar{F}_T =\projss\circ\phi_T\circ K_A =\phi_T\circ\projss\circ K_A,
\end{equation}
where $\projss(x) = y$ if and only if $x\in\varsub{W}{ss}(y)$, {\em
i.e.} $\projss:\Basin{\gamma}\to\gamma$ maps the basin of $\gamma$ onto
the limit cycle $\gamma$ along strong-stable manifolds.  The abstract
notations (and notions) have their uses: the projection $\projss$
encapsulates the properties of the strong-stable manifolds, {\em e.g.}
by definition, the strong stable manifold $\varsub{W}{ss}(\projss(x))$
passes through $x$ for any $x\in\Basin{\gamma}$, and the nullspace of
the Jacobian matrix $D\projss(x)$ is precisely the tangent space of the
strong stable manifold $\varsub{W}{ss}(\projss(x))$ at
$x$.\footnote{Note that the commutation relation $\phi_t\circ\projss
=\projss\circ\phi_t$ expresses the invariance of the strong-stable
foliation under $\phi_t$.  The map $K_A$, in general, has nothing to do
with the flow and does not commute with the other maps.}

\paragraph{Algorithm (Phase resetting curves via stable subspaces).}
\begin{enumerate}

\item Estimate the period $T_0$ of the limit cycle $\gamma$ by
  numerically solving the unforced equations starting with a point on or
  near $\gamma$.

\item\label{alg-osel:step1} Discretize $\gamma$ by subdividing the time
  interval $[0,T_0)$ into $N$ intervals and computing the
  corresponding points $x_i\in\gamma.$ Fix an arbitrary reference point
  $x_0$ on $\gamma$ so that each point on $\gamma$ can be assigned a
  unique phase $\theta\in[0,T_0)$.

\item For each point $x_i$ computed in the previous step, compute the
  Jacobian $DH(x_i)$ of the {\HH} flow field $H$ at that point.

\item\label{alg-osel:step4} Using the results of the previous two steps,
solve
\begin{equation}
\begin{array}{cl}
\dot{x} &= -H(x),\\
\dot{\xi} &= \eta - \iprod{\eta}{\xi}\xi,\\
\eta &= DH(x)^T\xi,
\end{array}
\end{equation}
using the grid points $\set{x_i}$ computed in the previous steps.  The
$\dot{x}$ part of the equation above is clearly numerically unstable,
but that is not a problem because we have already have a numerical
representation of $\gamma.$

The equations above are a variant of the usual method for computing
Lyapunov exponents {\cite{gh}}.  They preserve the length of $\xi(t)$,
though in practice, it is necessary to rescale $\xi(t)$ to ensure that
this constraint is maintained.  As $t\to\infty,$ $\xi(t)$ becomes
orthogonal to the strong-stable linear subspace $\varsub{E}{ss}(x(t))$
of $\gamma$.  The subspace $\varsub{E}{ss}(x(t))$ is tangent to the
strong-stable manifold $\varsub{W}{ss}(x(t))$ at $x(t)$.\footnote{These
equations can be generalized to the following:
\begin{equation}
\begin{array}{cl}
\dot{x} =& -H(x)\\
\dot{\xi}_i =& \eta_i -\iprod{\eta_i}{\xi_i}\xi_i -
               \sum_{j<i}\of{\iprod{\xi_i}{\dot{\xi}_j} + 
                             \iprod{\eta_i}{\xi_j}}\xi_j\\
\eta_i =& DH(x)^T\xi_i
\end{array}
\end{equation}
If the vectors $(\xi_i)$ form an orthonormal basis at $t=0$, then the
equations will guarantee that $(\xi_i(t))$ are orthonormal for all $t >
0$.  Again, it will be necessary to perform Gram-Schmidt
orthogonalizations periodically to maintain this constraint.  The vector
$\xi_1(t)$, as before, converges to a vector orthogonal to
$\varsub{E}{ss}(x(t))$.  So $\of{\xi_2(t),\xi_3(t),\xi_4(t)}$ span
$\varsub{E}{ss}(x(t)).$ Similarly, the vectors $(\xi_3(t),\xi_4(t))$
span the subspace consisting of the 2 fastest contracting directions,
and $(\xi_4(t))$ spans the fastest contracting direction.}

\item Using Eq.~\ref{eqn:phase-resetting-def-too} in combination with
  the linear subspaces computed in the previous step, we can now
  approximate the phase resetting curve.  Start with a point
  $x\in\gamma$ and compute $\Phi_t(K_A(x))$ for increasing $t$.  Let
  $t_0 > 0$ be the minimum positive time at which $\Phi_{t_0}(K_A(x))$
  has (i) returned to a small, fixed neighborhood of $\gamma$ (in this
  paper this is chosen to be a neighborhood of distance $10^{-4}$ around
  $\gamma$); and (ii) $\Phi_{t_0}(K_A(x))$ lies within one of the
  pre-computed linear subspaces $\varsub{E}{ss}(x_*)$ for some point
  $x_*$.  Let $\theta_*$ denote the phase of the point $x_*$.  Then the
  new phase of the system is $(T+\theta_*-t_0)\mbox{ (mod $T_0$)}$.

\item\label{alg-osel:step6} Proceed to the next grid point and repeat.

\item When the derivative of the phase resetting curve becomes large or
  infinite, it may be necessary to adaptively generate the grid points
  on which the curve is evaluated.  Generally speaking, the grid
  $\set{x_i}$ constructed in Step {\ref{alg-osel:step1}} need not equal
  the grid $\set{x'_i}$ on which the phase resetting curve is evaluated.
  In particular, the grid $\set{x'_i}$ can be adaptively chosen to
  ensure that $\abs{\hat{f}_a(x'_{i+1}) - \hat{f}_a(x'_i)}\leq\eps$,
  where $\hat{f}_a$ denotes the computed phase resetting curve and
  $\eps$ is a fixed number, in this paper usually $0.1$.  This adaptive
  mechanism provides a way to detect discontinuities in $f_T$.

\end{enumerate}


\section*{Acknowledgements}

It is also a pleasure to thank Lai-Sang Young, Eric Shea-Brown, Adi
Rangan, Charlie Peskin, John Rinzel, and Louis Tao for many helpful and
stimulating conversations.  I am grateful to the referees for their help
in improving the exposition and for pointing out relevant references.
This work is partially supported by an NSF Postdoctoral Fellowship.











\bibliographystyle{siam}
\bibliography{main}

\begin{thebibliography}{10}

\bibitem{aihara}
{\sc K.~Aihara and G.~Matsumoto}, {\em Chaotic oscillations and bifurcations in
  squid giant axons}, in Chaos, Nonlinear Science: Theory and Applications,
  Manchester University Press, 1986.

\bibitem{alexander}
{\sc J.~C. Alexander, E.~J. Doedel, and H.~G. Othmer}, {\em On the resonance
  structure in a forced excitable system}, SIAM Journal on Applied Mathematics,
  50 (1990), pp.~1373--1418.

\bibitem{best}
{\sc E.~N. Best}, {\em Null space in the {H}odgkin-{H}uxley equations: a
  critical test}, Biophysical Journal, 27 (1979), pp.~87--104.

\bibitem{brown}
{\sc E.~Brown, J.~Moehlis, and P.~Holmes}, {\em On the phase reduction and
  response dynamics of neural oscillator populations}, Neural Computation, 16
  (2004), pp.~673--715.

\bibitem{cronin}
{\sc J.~Cronin}, {\em Mathematical Aspects of Hodgkin-Huxley Neural Theory},
  Cambridge University Press, 1987.

\bibitem{doi1}
{\sc S.~Doi}, {\em A chaotic map with a flat segment can produce a
  noise-induced order}, Journal of Statistical Physics, 55 (1989),
  pp.~941--964.

\bibitem{ermentrout}
{\sc G.~B. Ermentrout}, {\em Simulating, Analyzing, and Animating Dynamical
  Systems: A Guide to XPPAUT for Researchers and Students}, vol.~14 of
  Software, Environments, and Tools, SIAM, 2002.

\bibitem{ek}
{\sc G.~B. Ermentrout and N.~Kopell}, {\em Multiple pulse interactions and
  averaging in systems of coupled neural oscillators}, Journal of Mathematical
  Biology, 29 (1991), pp.~195--217.

\bibitem{ermentrout1}
{\sc G.~B. Ermentrout, M.~Pascal, and B.~Gutkin}, {\em The effects of spike
  frequency adaptation and negative feedback on the synchronization of neural
  oscillators}, Neural Computation, 13 (2001), pp.~1285--1310.

\bibitem{geist}
{\sc K.~Geist, U.~Parlitz, and W.~Lauterborn}, {\em Comparison of different
  methods for computing {L}yapunov exponents}, Progress of Theoretical Physics,
  83 (1990), pp.~875--893.

\bibitem{glass}
{\sc L.~Glass and M.~C. Mackey}, {\em From Clocks to Chaos: The Rhythms of
  Life}, Princeton University Press, 1988.

\bibitem{guckenheimer}
{\sc J.~Guckenheimer}, {\em Isochrons and phaseless sets}, Journal of
  Theoretical Biology, 1 (1974), pp.~259--273.

\bibitem{gh}
{\sc J.~Guckenheimer and P.~Holmes}, {\em Nonlinear Oscillations, Dynamical
  Systems, and Bifurcations of Vector Fields}, Springer-Verlag, 1983.

\bibitem{guttman}
{\sc R.~Guttman, L.~Feldman, and E.~Jakobsson}, {\em Frequency entrainment of
  squid axon membrane}, Journal of Membrane Biology, 56 (1980), pp.~9--18.

\bibitem{hayashi3}
{\sc H.~Hayashi, S.~Ishizuka, and K.~Hirakawa}, {\em Chaotic response of the
  pacemaker neuron}, Journal of the Physical Society of Japan, 54 (1985),
  pp.~2337--2346.

\bibitem{hayashi4}
{\sc H.~Hayashi, S.~Ishizuka, M.~Ohta, and K.~Hirakawa}, {\em Chaotic behavior
  in the {\em onchidium} giant neuron under sinusoidal stimulation}, Physics
  Letters, 88A (1982), pp.~435--438.

\bibitem{hayashi5}
{\sc H.~Hayashi, M.~Nakao, and K.~Hirakawa}, {\em Entrained, harmonic,
  quasiperiodic and chaotic responses of the self-sustained oscillation of {\em
  nitella} to sinusoidal stimulation}, Journal of the Physical Society of
  Japan, 52 (1983), pp.~344--351.

\bibitem{hh}
{\sc A.~L. Hodgkin and A.~F. Huxley}, {\em A quantitative description of
  membrane current and its application to conduction and excitation in nerve},
  Journal of Physiology, 117 (1952), pp.~500--544.

\bibitem{holden}
{\sc A.~V. Holden}, {\em The response of excitable membrane models to a cyclic
  input}, Biology and Cybernetics,  (1975), pp.~1--7.

\bibitem{hubel}
{\sc D.~H. Hubel}, {\em Eye, Brain, and Vision}, Scientific American Library,
  1988.

\bibitem{hunter-milton}
{\sc J.~D. Hunter, J.~G. Milton, P.~J. Thomas, and J.~D. Cowan}, {\em Resonance
  effect for neural spike time reliability}, Journal of Neurophysiology, 80
  (1998).

\bibitem{ledrappier-young}
{\sc F.~Ledrappier and L.-S. Young}, {\em Entropy formula for random
  transformations}, Probability Theory and Related Fields, 80 (1988),
  pp.~217--240.

\bibitem{mainen}
{\sc Z.~F. Mainen and T.~J. Sejnowski}, {\em Reliability of spike timing in
  neocortical neurons}, Science, 268 (1995), pp.~1503--1506.

\bibitem{matsumoto1}
{\sc G.~Matsumoto, K.~Aihara, Y.~Hanyu, N.~Takahashi, S.~Yoshizawa, and J.~ichi
  Nagumo}, {\em Chaos and phase locking in normal squid axons}, Physics Letters
  A, 123 (1987), pp.~162--166.

\bibitem{matsumoto2}
{\sc G.~Matsumoto, N.~Takahashi, and Y.~Hanyu}, {\em Chaos, phase locking and
  bifurcation in normal squid axons}, in Chaos in Biological Systems, Plenum,
  1987.

\bibitem{nakao}
{\sc H.~Nakao, K.~suke Arai, K.~Nagai, Y.~Tsubo, and Y.~Kuramoto}, {\em
  Synchrony of limit-cycle oscillators induced by random external impulses},
  Physical Review E, 72 (2005).

\bibitem{wang}
{\sc A.~Oksasoglu and Q.~Wang}, {\em Strange attractors in periodically-kicked
  {C}hua's circuit}, International Journal of Bifurcation and Chaos,  (to
  appear).

\bibitem{oprisan-canavier}
{\sc S.~A. Oprisan and C.~C. Canavier}, {\em The influence of limit cycle
  topology on the phase resetting curve}, Neural Computation, 14 (2002),
  pp.~1027--1057.

\bibitem{pakdaman}
{\sc K.~Pakdaman and S.~Tanabe}, {\em Random dynamics of the {H}odgkin-{H}uxley
  neuron model}, Physical Review E, 64 (2001).

\bibitem{recipes}
{\sc W.~H. Press, B.~P. Flannery, S.~A. Teukolsky, and W.~T. Vetterling}, {\em
  Numerical Recipes in C}, Cambridge University Press, 1992.

\bibitem{rieke}
{\sc F.~Rieke}, {\em Spikes: Exploring the Neural Code}, MIT Press, {\em
  c}1997.

\bibitem{takabe}
{\sc T.~Takabe, K.~Aihara, and G.~Matsumoto}, {\em Response characteristics of
  the {H}odgkin-{H}uxley equations to pulse-train stimulation}, Transactions of
  IEICE Japan (Japanese edition), J71-A (1988), pp.~744--750.

\bibitem{wy1}
{\sc Q.~Wang and L.-S. Young}, {\em Strange attractors with one direction of
  instability}, Communications in Mathematical Physics, 218 (2001).

\bibitem{wy2}
\leavevmode\vrule height 2pt depth -1.6pt width 23pt, {\em From invariant
  curves to strange attractors}, Communications in Mathematical Physics, 225
  (2002).

\bibitem{wy3}
\leavevmode\vrule height 2pt depth -1.6pt width 23pt, {\em Strange attractors
  in periodically-kicked limit cycles and {H}opf bifurcations}, Communications
  in Mathematical Physics, 240 (2003).

\bibitem{wy4}
\leavevmode\vrule height 2pt depth -1.6pt width 23pt, {\em Toward a theory of
  rank-1 attractors}, preprint,  (2004).

\bibitem{williams-bowtell}
{\sc T.~L. Williams and G.~Bowtell}, {\em The calculation of frequency-shift
  functions for chains of coupled oscillators, with application to a network
  model of the lamprey locomotor pattern generator}, Journal of Computational
  Neuroscience, 4 (1997), pp.~47--55.

\bibitem{winfree}
{\sc A.~Winfree}, {\em The Geometry of Biological Time, Second Edition},
  Springer-Verlag, 2000.

\bibitem{young}
{\sc L.-S. Young}, {\em Ergodic theory of differentiable dynamical systems}, in
  Real and Complex Dynamics, NATO ASI Series, Kluwer Academic Publishers, 1995,
  pp.~293--336.

\bibitem{zaslavsky}
{\sc G.~Zaslavsky}, {\em The simplest case of a strange attractor}, Physics
  Letters, 69A (1978), pp.~145--147.

\bibitem{zhou-kurths}
{\sc C.~Zhou and J.~Kurths}, {\em Noise-induced synchronization and coherence
  resonance of a {H}odgkin-{H}uxley model of thermally sensitive neurons},
  Chaos, 13 (2003), pp.~401--409.

\end{thebibliography}
\end{document}